\documentclass[12pt,a4paper,english]{article}
\usepackage{amsmath}
\usepackage{amssymb}
\usepackage{amsthm}
\usepackage{bbm}
\usepackage[authoryear]{natbib}
\usepackage{booktabs}
\usepackage{graphicx}
\usepackage{epstopdf}
\usepackage{caption}
\usepackage{subcaption}

\usepackage{babel}
\usepackage{color}
\usepackage{enumitem}
\usepackage{rotating}
\usepackage{pdflscape}
\usepackage{hyphenat}
\usepackage[a4paper,margin=1.25in]{geometry}
\usepackage{algorithm}
\usepackage{algpseudocode}
\usepackage{float}
\usepackage[]{threeparttable}
\usepackage{multirow}
\newcommand\fnote[1]{\captionsetup{font=footnotesize}\caption*{#1}}
\usepackage[onehalfspacing,nodisplayskipstretch]{setspace}

\usepackage{xcolor}
\usepackage[title]{appendix}
\definecolor{lightblue}{rgb}{0, 0.4, 0.6}

\usepackage[capitalise,nameinlink,noabbrev]{cleveref}
\crefname{equation}{}{}
\crefname{enumi}{}{}
\crefname{appsec}{Appendix}{Appendices}
\crefname{appsubsec}{Appendix}{Appendices}

\newtheorem{nondecr}{Theorem}
\newtheorem{theorem}[nondecr]{Theorem}
\newtheorem{nondecrl}{Lemma}
\newtheorem{lemma}[nondecrl]{Lemma}
\newtheorem{nondecrd}{Definition}
\newtheorem{definition}[nondecrd]{Definition}

\newtheorem{manualtheoreminner}{Assumption}
\newenvironment{manualassumption}[1]{%
  \renewcommand\themanualtheoreminner{#1}%
  \manualtheoreminner
}{\endmanualtheoreminner}

\begin{document}

\title{Scalable simultaneous inference in high-dimensional linear regression models}
\author{Tom Boot
\thanks{Department of Economics, Econometrics and Finance, University of Groningen,  Nettelbosje 2,
9747 AE Groningen,
The Netherlands, e-mail: \textsf{t.boot@rug.nl}}\\
University of Groningen
	 \and Didier Nibbering\thanks{Department of Econometrics \& Business Statistics, Monash University, Clayton VIC 3800, Australia, e-mail: \textsf{didier.nibbering@monash.edu} \newline
	\indent
	We would like to thank Paul Bekker, Patrick Groenen, Trevor Hastie, Christiaan Heij, Richard Paap, Andreas Pick, and participants of the workshop on statistical learning and econometrics at the Erasmus University Rotterdam for helpful comments. }\\
Monash University}

\maketitle
\begin{abstract}
The computational complexity of simultaneous inference methods in high-dimensional linear regression models quickly increases with the number variables. This paper proposes a computationally efficient method based on the Moore-Penrose pseudoinverse. Under a symmetry assumption on the available regressors, the estimators are normally distributed and accompanied by a closed-form expression for the standard errors that is free of tuning parameters. We study the numerical performance in Monte Carlo experiments that mimic the size of modern applications for which existing methods are computationally infeasible. We find close to nominal coverage, even in settings where the imposed symmetry assumption does not hold. Regularization of the pseudoinverse via a ridge adjustment is shown to yield possible efficiency gains.
\end{abstract}
\newpage

\section{Introduction}\label{sec: introduction}
Modern data sets require inferential procedures that are valid in high-dimensional settings where the number of variables $p$ exceeds the sample size $n$. Typical examples occur in genome wide association studies, such as \citet{candes2018knockoffs} where $n=O(10^3)$ and $p=O(10^5)$, or \citet{rietveld2013gwas} where $n=O(10^5)$ and $p=O(10^6)$. 

Debiasing, or desparsification, strategies allow simultaneous inference on full parameters vector in high-dimensional linear regression models.
Debiasing relies on an approximate inverse of the singular high-dimensional covariance matrix of the regressors, and the use of the lasso to correct for the bias resulting from this approximation. The remaining bias, generally factorized in a part due to the approximate inverse and one due to the lasso estimation error, is then negligible compared to the variance of the estimator. 

One downside of current debiasing strategies is that the computational burden increases rapidly with the number of included variables, as the approximate inverse is constructed by solving a series of $p$ optimization problems. \citet{javanmard2014confidence} rely on direct numerical optimization for each variable to find an approximate inverse, while \citet{van2014asymptotically,zhang2014confidence,ning2017general,dezeure2017high} solve $p$ lasso problems to estimate the inverse covariance matrix.
Even with modern computing power, this puts a strain on computational resources when $p$ is large. As a result, existing simulation studies are typically limited to settings where $(n,p)=O(10^2)$. This drives a wedge between the scale demanded by modern applications and the scale at which the methods can be tested.

This paper proposes a scalable debiasing method that allows for simultaneous inference in high-dimensional linear regression models. The required approximate inverse is constructed by a scaled Moore-Penrose pseudoinverse, where the scaling ensures that the approximation bias becomes negligible. Replacing the sequence of optimization problems with a single pseudoinverse lowers the computational complexity with several orders of magnitude in the number of variables. Moreover, standard errors are available in closed form and free of tuning parameters.

The confidence intervals of the Moore-Penrose pseudoinverse estimator are shown to be valid under a symmetry assumption on the variables, which for example allows variables to be independent, equicorrelated, or to follow certain factor structures. We show that a ridge adjustment of the Moore-Penrose pseudoinverse estimator can lead to a power gain. 

We carry out a simulation study where we consider dimensions $n$ and $p$ for which existing methods have substantial computational costs. We find that coverage is close to nominal for settings where the variables satisfy the symmetry assumption. The proposed ridge adjustment ensures that the method is competitive to existing alternatives in terms of power. We find undercoverage on the set of nonzero coefficients for some, but not all, of the settings that violate the symmetry assumption.

The proposed estimator is related to the corrected ridge estimator of \citet{buhlmann2013statistical}, who finds that under a fixed design
an additional bias correction is required that leads to conservative inference. Our estimator differs in 
the scaling factor that renders the approximation bias due to the diagonal elements of the approximate inverse exactly equal to zero. Under a random design and a slightly stronger sparsity assumption, we show that the additional bias correction is no longer required.



The outline of the paper is as follows. \cref{sec: gf} introduces the estimation approach and the proposed estimators. The theoretical properties of the Moore-Penrose pseudoinverse and ridge regression are presented in \cref{sec: theory}. \cref{sec: mc} illustrates these results through Monte Carlo simulations. \cref{sec: conclusion} concludes.

\textbf{Notation}$\quad$ We use the following notation throughout the paper: For any $n\times 1$ vector $a=(a_1,\dots,a_n)^\top$, the $l_{q}$-norm is defined as $||a||_q:=(\sum^n_{i=1} |a_i|^q)^{1/q}$ for $q>0$ and $||a||_{0}$ denotes the number of nonzero elements of $a$. The maximum norm is written as $||a||_{\infty} = \max(|a_1|,\dots,|a_n|)$. For a $p\times n$ matrix $A$, the $l_{q}$-norm is defined as $||A||_q:=\sup_{x, ||x||_{q}=1}\left\{||Ax||_{q}\right\}$ and the maximum norm is written as $||A||_{\max} = \max_{i=1,\dots,n,j=1,\dots,p}|A_{ij}|$. When $A$ is $n\times n$, we write $\lambda_{\min}(A)$ and $\lambda_{\max}(A)$ for the minimum and maximum eigenvalues of $A$. The $n\times n$ identity matrix is denoted by $I_{n}$. The vector $e_{i}$ has its $i$-th entry equal to 1 and zeros everywhere else. For the regressor matrix $X$, we index the rows with the subscript $i=1,\ldots,n$ and the columns with the subscript $j=1,\ldots,p$. If $U$ is a $p\times p$ orthogonal matrix, we write $U\in \mathcal{O}(p)$. When two random variables $X$ and $Y$ follow the same distribution, this is denoted as $X\overset{(d)}{=}Y$.

\section{High-dimensional linear regression}\label{sec: gf}
Consider the high-dimensional linear model
\begin{equation}\label{eq:dgp}
y = X\beta + \varepsilon, \quad \varepsilon \sim N(0,\sigma^2 I_{n}),
\end{equation}
where $y$ is an $n \times 1$ response vector, $X$ an $n \times p$ regressor matrix, $\beta=(\beta_1,\dots,\beta_p)^\top$ a $p \times 1 $ vector of unknown regressor coefficients, and $\varepsilon$ an $n \times 1$ vector of errors which are $N(0,\sigma^2 I_n)$ and independent of $X$.

We develop a computationally efficient debiased estimator $\hat{\beta}$ for $\beta$ when $p>n$, with a closed-form expression for the covariance matrix of the estimator. Standard errors and confidence intervals for $\hat{\beta}$ are free of tuning parameters, closed-form, and scalable to a very large number of variables. Below we develop the estimator. \cref{sec:inference} provides an algorithmic overview of the proposed inference method.

\subsection{A debiased estimator}
We consider estimators for $\beta$ of the form
\begin{equation}\label{eq:start}
\begin{split}
    \hat{\beta}^{b} &= \hat{M}y\\
    &=\hat{M}X\beta +\hat{M}\varepsilon\\
    & = \beta +\left(\hat{M}X - I_{p}\right)\beta + \hat{M}\varepsilon,\\
\end{split}
\end{equation}
where $\hat{M}$ is a $p\times n$ matrix. The second term of \eqref{eq:start} represents a bias which depends on the choice of $\hat{M}$. When $p\leq n$ and $\text{E}[\varepsilon|X]=0$, ordinary least squares yields unbiased estimators by choosing $\hat{M} = (X^\top X)^{-1}X^{\top}$.
When $p>n$, the matrix $X^\top X$ is singular, and any choice of $\hat{M}$ will induce bias.

Given an initial estimator $\hat{\beta}^{\text{init}}$, we can reduce the bias in \eqref{eq:start} by applying the correction
\begin{equation}\label{eq:betac}
\begin{split}
    \hat{\beta} &= \hat{\beta}^{b} - \left(\hat{M}X - I_{p}\right)\hat{\beta}^{\text{init}} \\
    & = \beta +\left(\hat{M}X - I_{p}\right)\left(\beta-\hat{\beta}^{\text{init}}\right) + \hat{M}\varepsilon.\\
\end{split}
\end{equation}

For the initial estimator $\hat{\beta}^{\text{init}}$ we use the lasso estimator of \citet{tibshirani1996regression} defined by
\begin{equation}\label{eq:lasso}
\hat{\beta}^{\text{lasso}} = \arg\min_{b}\left[\frac{1}{n}(y-Xb)^\top(y-Xb) + \lambda||b||_{1}\right].
\end{equation}
This leads to an interpretation of $\hat{\beta}$ as a `desparsified' version of the lasso estimator \citep{van2014asymptotically}, since we can write $
    \hat{\beta} = \hat{\beta}^{\text{lasso}}+\hat{M}(y-X\hat{\beta}^{\text{lasso}})$.

\subsection{The approximate inverse $\hat{M}$}\label{sec: chooseM}
The goal of this paper is to introduce choices of $\hat{M}$ that are computationally efficient and for which the bias of the corrected estimator is small. For the bias to be small, $\hat{M}X$ has to be close to the $p\times p$ identity matrix $I_{p}$. Hence, we refer to $\hat{M}$ as an approximate inverse for $X$.

We ensure that the diagonal terms of $\hat{M}X - I_{p}$ are identically equal to zero by introducing a $p\times p$ diagonal matrix $\hat{D}$, with diagonal elements $\hat{d}_{j}$, and taking
\begin{equation}\label{eq:diagonalscaling}
\hat{M} = \hat{D}\tilde{M}, \qquad
\hat{D}_{jj}  = (e_{j}^\top\tilde{M}Xe_{j})^{-1}.
\end{equation}
The effectiveness of $\hat{M}$ as an approximate inverse is then determined by the magnitude of the off-diagonal elements of $\hat{M}X$.

\subsubsection{Moore-Penrose pseudoinverse}
As a tuning parameter free choice for $\tilde{M}$ in \eqref{eq:diagonalscaling} we consider the Moore-Penrose pseudoinverse (MPI) given by $X^\top(XX^\top)^{-1}$. That is,
\begin{align}\label{eq: mppi}
   \hat{M} = \hat{D}X^\top(XX^\top)^{-1}.
\end{align}
The diagonal elements $\hat{D}_{jj}$ of the diagonal scaling matrix $\hat{D}$ equal
\begin{equation}\label{eq: diagmppi}
\hat{D}_{jj}= \left[e_{j}^\top X^\top(X X^\top)^{-1}Xe_{j}\right]^{-1}.
\end{equation}
This provides a closed-form expression for the approximate inverse. Since the bias term of the estimator is of lower order compared to the variance, the covariance of $\hat{\beta}$ is available in closed form as well,
\begin{equation}\label{eq:varbeta}
V(\hat{\beta}) =\sigma^2 \hat{D}X^\top(XX^\top)^{-2}X\hat{D}.
\end{equation}

\subsubsection{Ridge regularization}
We also consider the use of a ridge adjustment (RID), which in finite samples can lead to a power gain. Define
\begin{align}\label{eq: ridge}
    \hat{M} = \hat{D}(X^\top X+\gamma I_{p})^{-1}X^\top,
\end{align}
where $\gamma$ denotes the ridge penalty and the elements of the diagonal scaling matrix $D$ equal
\begin{equation}\label{eq: diagridge}
\hat{D}_{jj}=\left(e_{j}^\top(X^\top X+\gamma I_{p})^{-1}X^\top Xe_{j}\right)^{-1}.
\end{equation}

The ridge-adjustment is the natural regularization procedure as we can define the Moore-Penrose pseudoinverse as
\begin{equation}\label{eq:ridgeMP}
\begin{split}
X^\top(XX^\top)^{-1} &= \lim_{\gamma\rightarrow 0}\left(X^\top X + \gamma I_{p}\right)^{-1}X^\top\\
& = \lim_{\gamma\rightarrow 0}X^\top\left(XX^\top + \gamma I_{n}\right)^{-1},
\end{split}
\end{equation}
see e.g. \citet{albert1972regression}.

\subsubsection{Alternative specifications}
\citet{zhang2014confidence} develop the estimator $\hat{\beta}$ in \eqref{eq:betac} with $\hat{M}=\frac{1}{n}\bar{M}X^\top$, where $\bar{M}$ is an approximate inverse to the empirical covariance matrix $\frac{1}{n} X^\top X$. This inverse is found by a series of lasso regressions. \citet{van2014asymptotically} generalize this method to deal with nonlinear models and prove semiparametric optimality. \citet{javanmard2014confidence} use direct numerical optimization to obtain $\bar{M}$.

The ridge adjusted estimator, but with a different scaling factor, is studied by \citet{buhlmann2013statistical} for a fixed design and a less restrictive sparsity assumption then considered in \cref{sec: theory}. In this case an additional bias correction is necessary that results in conservative inference.

\subsection{Inference}\label{sec:inference}
With a specification for the approximate inverse $\hat{M}$, confidence intervals can be constructed for the $j$-th element of $\beta$ in \eqref{eq:dgp} as
\begin{equation}\label{eq: ci}
\left[\hat{\beta}_{j} - z_{\alpha/2}\sqrt{\hat{\Omega}_{jj}/n},\quad\hat{\beta}_{j} + z_{\alpha/2}\sqrt{\hat{\Omega}_{jj}/n}\right],
\end{equation}
where $\hat{\beta}_{j}$ is the $j$-th element of $\hat{\beta}$ in \eqref{eq:betac}, $z_{\alpha/2}$ is the $\alpha/2$ critical value for the standard normal distribution, and $\sqrt{\hat{\Omega}_{jj}/n}$ is the standard error of $\hat{\beta}_{j}$.

The covariance matrix of the estimator $\hat{\beta}$ is equal to $\hat{\Omega}=\hat{\sigma}^2 n\hat{M}\hat{M}^\top$. The noise level $\sigma^2$ is estimated with the scaled lasso estimator of \citet{scaledlassoSun}, defined as
\begin{equation}\label{eq:scaledlasso}
(\hat{\beta}^{\text{scaled}},\hat{\sigma}^{\text{scaled}}) = \arg\min_{b,s}\left[\frac{1}{2ns}(y-Xb)^\top(y-Xb) + \frac{s}{2}+\lambda_0||b||_{1}\right],
\end{equation}
with $\lambda_0=\sqrt{2\log(p)/n}$. This scaled lasso estimator is commonly used for noise estimation in high-dimensional regression settings.

\cref{alg:method} outlines the proposed inference method. The algorithm requires a value for the ridge penalty. The numerical results in \cref{sec: mc} indicate that the ridge adjustment only has an effect when $n$ and $p$ are of the same order of magnitude. In this case, following our theoretical results, setting  $\gamma = p\sqrt{\log p /n}$ improves the accuracy compared to the MPI estimator, while keeping coverage at the nominal level.

\begin{algorithm}[t]
\caption{Inference in high-dimensional linear regression models}\label{alg:method}
\begin{algorithmic}[1]
\Require Data $(y,X)$, significance level $\alpha$, and regularization strength $\gamma$.
\State Set $z_{\alpha/2}$ to the $\alpha/2$ critical value for the standard normal distribution.
\State Get $\hat{\beta}^{\text{lasso}}$ from \eqref{eq:lasso} with $\lambda$ selected through 100-fold cross-validation.
\State Get $\hat{\sigma}^{\text{scaled}}$ from \eqref{eq:scaledlasso}.
\If {$\gamma=0$}
\State Set $\hat{D}_{jj}= \left[e_{j}^\top X^\top (X X^\top )^{-1}Xe_{j}\right]^{-1}$ and $\hat{D}_{ij}=0$ if $i\neq j$
\State   $\hat{M}=\hat{D}X^\top (XX^\top )^{-1}$
\Else
\State Set $\hat{D}_{jj}=\left(e_{j}^\top (X^\top X + \gamma I_{p})^{-1}X^\top Xe_{j}\right)^{-1}$ and $\hat{D}_{ij}=0$ if $i\neq j$
\State  $\hat{M}=\hat{D}(X^\top X+\gamma I_{p})^{-1}X^\top $
\EndIf
\State $\hat{\beta}=\hat{M}y - \left(\hat{M}X - I_{p}\right)\hat{\beta}^{\text{lasso}}$
\State $\hat{\Omega}=(\hat{\sigma}^{\text{scaled}})^2n\hat{M}\hat{M}^\top $
\State \Return CI $\left[\hat{\beta}_{j}^{c} - z_{\alpha/2}\sqrt{\hat{\Omega}_{jj}/n},\quad\hat{\beta}_{j}^{c} + z_{\alpha/2}\sqrt{\hat{\Omega}_{jj}/n}\right] \forall j$
\end{algorithmic}
\end{algorithm}

\cref{alg:method} uses the lasso as the initial estimator and the scaled lasso for noise estimation. However, alternative initial estimators, such as proposed by \citet{caner2014asymptotically}, can also be used as long as they satisfy a sufficiently tight accuracy bound on  $||\beta-\hat{\beta}^{\text{init}}||_{1}$. Similarly, any consistent estimator for $\sigma$ can be used for noise estimation. Instead of the scaled lasso, \citet{reid2013study} estimate the error variance from the lasso residuals with a degree of freedom correction. \citet{yu2019estimating} take the minimizing value of the lasso regression problem.

\subsection{Computational complexity}
\cref{tab:ccc} shows that the computational complexity of the proposed methods is several orders of magnitude in the number of variables lower compared to existing methods.

  \begin{table}[tb!]
 	\caption{\label{tab:ccc} Computational complexity comparison\vspace{0.1cm}}
	\begin{tabular}{llll}
		\toprule
		Method & Complexity& Method & Complexity\\
		\midrule
		MPI   &  $n \cdot p^2  $ &\citet{van2014asymptotically} & $p^4$ \\
		RID  &  $n \cdot p^2 $&\citet{javanmard2014confidence} & $p^4$\\
		\bottomrule
	\end{tabular}
\end{table}

For the Moore-Penrose pseudoinverse, the leading order term in the computation complexity is the matrix multiplication $X(XX^\top )^{-1}X^\top $, which costs $O(n\cdot p^2)$. The same holds for the ridge regularized version.

The method of \citet{van2014asymptotically} solves a penalized regression problem with $p-1$ explanatory variables for every column of the empirical covariance matrix. A fast solver is provided by the lars algorithm of \citet{efron2004least}, which has a complexity of $O(p^3)$ for each column. Therefore, the complexity for obtaining the approximate inverse covariance matrix equals $O(p^4)$ in total. \citet{javanmard2014confidence} state that this is equivalent to the complexity of their method.

All methods use an initial estimator, which is provided by the lasso. With $r = \min\left\{n,p\right\}$, this costs $O(r^3)$ at a fixed value of the penalization parameter. For the initial estimator we apply $K$-fold cross-validation, and $r = n$, such that the computational complexity of all methods is at least $O(K\cdot n^3)$. However, for none of the methods this is the leading order term.

\section{Theoretical results}\label{sec: theory}
\subsection{Assumptions}
We make the following assumptions on the expansion rates of $n$ and $p$, the regressor matrix $X$ and the errors $\varepsilon$.

\begin{manualassumption}{1}\label{ass:sparsity} The dimensions $n$ and $p$ are such that $C_{1}\log p<n<C_{2}p$ for some sufficiently large constant $C_{1}>0$ and the constant $C_{2}<1$. The sparsity $s_{0} \equiv ||\beta||_{0}$ satisfies $s_{0} = o\left(\frac{\sqrt{n}}{\log p}\right)$.
\end{manualassumption}

\begin{manualassumption}{2a}\label{ass:R}
Let $R$ be an $n\times p$ matrix with independent subgaussian rows $r_{i} = (\Sigma^{R})^{1/2}\tilde{r}_{i}$, where  $\tilde{r}_{i}$ has independent subgaussian elements. The matrix $\Sigma^{R}=\emph{E}[R R^{\top}]$ satisfies $\max_{i=1,\ldots,p}\Sigma^{R}_{ii}=O(1)$, and $1/\lambda_{\min}(\Sigma^{R})=O(1)$. Assume $X=R$.
\end{manualassumption}

\begin{manualassumption}{2b}\label{ass:factor}
Let $F$ be an $n\times k$ matrix with independent rows and identically subgaussian distributed columns. Let $\Lambda$ be an $p\times k$ matrix with identically subgaussian distributed columns, and $\Sigma^{\Lambda} = \emph{E}[\Lambda\Lambda^{\top}]$ has finite eigenvalues bounded away from zero by a constant. Assume the factor structure $X=F\Lambda^{\top} + R$, with $k$ fixed, the matrices $R$, $F$ and $\Lambda$ independent, and $R$ as defined in \cref{ass:R} with the additional assumption that $\lambda_{\max}(\Sigma^{R})<M$.
\end{manualassumption}

\begin{manualassumption}{3}\label{ass:symmetry}
The rows of $X$ are spherically symmetric: The rows of $R$ satisfy $r_{i}\overset{(d)}{=}Tr_{i}$, and in case of a factor structure the columns of $\Lambda$ also satisfy $\lambda_{i}\overset{(d)}{=}T\lambda_{i}$, for any fixed $T\in \mathcal{O}(p)$.
\end{manualassumption}


\begin{manualassumption}{4}\label{ass:error}
The errors $\varepsilon\sim N(0,\sigma^2 I_{n})$ and $\varepsilon$ is independent of $X$.
\end{manualassumption}

\cref{ass:sparsity} imposes a condition on the sample size relative to the number of variables. We  operate under the assumption that $p>n$, but we assume that $n$ is sufficiently larger than $\log p$. The sparsity constraint restricts the number of non-zero coefficients in $\beta$ by $s_{0}=||\beta||_{0}$. The assumption is standard for inference in high-dimensional regression models \citep{van2014asymptotically,javanmard2014confidence}, but slightly stronger than the $s_{0}^2 = o\left(n/\log p\right)$ required for lasso consistency \citet{van2008high}. This discrepancy has been studied in more detail in \citet{javanmard2018debiasing}. Since our results only depend on the $l_{1}$ norm of the lasso estimation error, the sparsity assumption can be relaxed to allow for approximate sparsity \citep{chernozhukov2015valid}.

\cref{ass:R} is taken from Theorem 2.4 in \citet{javanmard2014confidence}. This assumption is necessary to bound the bias of the initial estimator. It allows for a random $X$ with subgaussian independent rows, but the subgaussian assumption rules out a factor structure for $X$, which is allowed by \cref{ass:factor}.

\cref{ass:symmetry} is necessary to bound the bias of the approximate inverse. It is for example satisfied when the columns of $X$ have an equicorrelated covariance matrix, but also covers classes of factor models. \citet{fan2008sure}
use this assumption to derive accuracy bounds on the Moore-Penrose pseudoinverse in high-dimensional variable screening. These results depend on the behavior of the right singular vectors under spherical symmetry. In a slightly different context, \citet{shah2020right} show that under spherical symmetry, the right singular vectors can be used to estimate high-dimensional covariance matrices that contain a factor component. We find in \cref{sec: mc} that our methods continue to work well in a range of settings in which the spherical symmetry assumption is violated, for example when the columns of $X$ have a covariance matrix with a Toeplitz structure.

Finally, \cref{ass:error} ensures that the estimated parameters are normally distributed. This assumption can be relaxed as shown by \citet{van2014asymptotically}, in which case valid inference can be achieved for a fixed number of coefficients.

\subsection{Theorems}
The following theorem presents the main result of this paper.
\begin{theorem}\label{theorem:main}
Suppose \cref{ass:sparsity} -- \ref{ass:error} hold. Let $\hat{\beta}=\hat{M}y - \left(\hat{M}X-I_p\right)\hat{\beta}^{\text{lasso}}$, with $\hat{M}$ as defined in \eqref{eq: mppi} and $\hat{D}$ as in \eqref{eq: diagmppi}.
Then,
\begin{align*}
    \sqrt{n}(\hat{\beta}-\beta) & = Z + o_{p}(1),\\
    Z|X &\sim N\left(0, \hat{\Omega}\right),
\end{align*}
where $\hat{\Omega}= \hat{\sigma}^2n\hat{M}\hat{M}^\top $ and $\hat{\Omega}_{jj} = O_{p}(1)$.
\end{theorem}
The proof is given in \cref{app:prooftheoremmain}. Theorem~\ref{theorem:main} shows that the estimator $\hat{\beta}$ in \eqref{eq:betac} has a lower order bias compared to the variance. Standard errors decrease at the usual $n^{-1/2}$ rate. Theorem~\ref{theorem:main} allows for the construction of confidence intervals that are uniformly valid when $p$ is finite. Uniformity is guaranteed since the bound on the lasso estimator given in \cref{lem:acclasso} holds uniformly over all sets $S_{0}$ of size $s_{0} = o(\sqrt{n}/\log p)$, see \citet{van2014asymptotically} for a discussion.

Since the resulting covariance matrix of the estimator is available in closed form, efficient multiple testing procedures as in \citet{buhlmann2013statistical} can be employed, together with joint tests on estimated coefficients.

\cref{theorem:main2} shows that for an appropriate choice of the regularization parameter, the ridge adjusted estimator can also be used for inference.

\begin{theorem}\label{theorem:main2}
Suppose \cref{ass:sparsity} -- \ref{ass:error} hold. Let $\hat{\beta}(\gamma)=\hat{M}y - \left(\hat{M}X-I_p\right)\hat{\beta}^{\text{lasso}}$, with $\hat{M}$ as in \eqref{eq: ridge} and $\hat{D}$ as in \eqref{eq: diagridge}. Suppose $\gamma=O(p\sqrt{\log p/n})$.
Then,
\begin{align*}
    \sqrt{n}(\hat{\beta}(\gamma)-\beta) & = Z + o_{p}(1),\\
    Z|X &\sim N\left(0,\hat{\Omega}\right),
\end{align*}
where $\hat{\Omega}= \hat{\sigma}^2 n\hat{M}\hat{M}^\top $ and $\hat{\Omega}_{jj} = O_{p}(1)$.
\end{theorem}
The proof is given in \cref{app:prooftheoremmain2}.

The reason one would opt for the regularized variants despite the additional tuning parameter $\gamma$ is provided by the following theorem.
\begin{theorem}\label{theorem:power} Suppose \cref{ass:sparsity} -- \ref{ass:error} hold.  Denote the variance of the estimator $\hat{\beta}_{j}$ under a diagonal scaling matrix $\hat{D}$ by $\hat{\Omega}_{jj}^{\text{MPI}}(\hat{D})$ when $\hat{M}$ is as in \eqref{eq: mppi} and $\hat{\Omega}_{jj}^{\text{RID}}(\hat{D})$ when $\hat{M}$ is as in \eqref{eq: ridge}. For the choice of $\gamma$ as in \cref{theorem:main2}, we have
	\begin{equation}\label{eq:th4}
\hat{\Omega}_{jj}^{\text{RID}}(\hat{D})-\hat{\Omega}^{\text{MPI}}_{jj}(\hat{D})< 0.
\end{equation}
\end{theorem}
The proof is given in \cref{app:prooftheorempower}. The proof shows that the difference between variances of the RID and MPI estimators increases in $\gamma$. Since the proof of \cref{theorem:main2} in \cref{app:prooftheoremmain2} suggests that $\gamma$ can attain a higher value when $p$ is close to $n$, we expect a larger potential efficiency gain of ridge regularization in these settings.



Note that \cref{theorem:power} requires the regularized estimator and the estimator based on the Moore-Penrose pseudoinverse to use the same diagonal scaling matrix. Using $\hat{D}$ from \eqref{eq: diagmppi} for the Moore-Penrose inverse and $\hat{D}$ from \eqref{eq: diagridge} for the ridge regularized inverse, does not yield an ordering in terms of power. However, in all cases we have encountered, the inequality in \cref{theorem:power} is satisfied when using the diagonal matrix specific to the estimator under consideration. This is also evident from the Monte Carlo results in \cref{sec: mc}.

\section{Monte Carlo Experiments}\label{sec: mc}
This section examines the finite sample behaviour of the proposed estimators in a Monte Carlo experiment. The performance of the estimators is compared with \citet{van2014asymptotically} (GBRD) and \citet{javanmard2014confidence} (JM), two existing methods for constructing confidence intervals in high-dimensional regression for all coefficients. All estimators use the lasso as initial estimator with a penalty term that minimizes the mean squared error under tenfold cross-validation. The series of lasso regressions in GBRD also use tenfold cross-validation, and we set the tuning parameter $\mu=2 \sqrt{n^{-1} \log p}$ in JM, which is equal to the value used in their simulation studies. 

The estimators are evaluated on mean absolute bias (MAE), coverage rate, power, and family-wise error rate (FWER). For each Monte Carlo replication, we calculate 
\begin{equation}
    \begin{split}
        \text{MAE} & = \frac{1}{s_{0}}\sum_{j\in S}|\hat{\beta}_{j}-\beta_{j}|,\\
        \text{Coverage} & = \frac{1}{s_{0}}\sum_{j\in S}\mathbbm{1}[{\beta}_{j}\in (\hat{\beta}_{j}-1.96 \sqrt{n^{-1}\hat{\Omega}_{jj}},\hat{\beta}_{j}+1.96 \sqrt{n^{-1}\hat{\Omega}_{jj}})],\\
        \text{Power}& = \frac{1}{s_{0}}\sum_{j\in S}\mathbbm{1}[H_{0,j} \text{ is rejected}],\\
        \text{FWER} &= \frac{1}{p-s_{0}}\sum_{j\in S^{c}} \mathbbm{1}[H_{0,j} \text{ is rejected}],
     \end{split}
\end{equation}
where MAE and coverage are similarly defined on the set of zero coefficients $S^{c}$, and $H_{0,j}\colon \beta_{j}=0$ for $j=1,\ldots,p$. We report the average MAE, coverage and power over the Monte Carlo replications. 

\subsection{Monte Carlo design}
The data generating process takes the form
\begin{equation}\label{eq: dgpmc}
y = X\beta + \varepsilon, \quad \varepsilon \sim N(0,\sigma^2 I_n),
\end{equation}	
where $X$ is an $n \times p$ regressor matrix, and $\beta$ a $p \times 1$ vector of coefficients. \cref{tab:design} shows the 14 simulation designs with different specifications for $\beta$ and $X$, as used in \citet{wang2015high}.

For each simulation design, 1000 data sets are generated with the number of predictors and sample size $(p,n)$ equal to $(200,100)$, $(1000,200)$, and $(10000,400)$. The error variance $\sigma^2$ is set to the value that satisfies $R^2=\text{var}(X\beta)/(\text{var}(X\beta)+\sigma^2)=50\%$. Note that due to the computational costs of GBRD and JAM, we only include MPI and RID in the Monte Carlo experiments with $(p,n)=(10000,400)$.

\begin{table}[tb!]
\centering\footnotesize
\caption{Simulation designs in the Monte Carlo experiments}\label{tab:design}
\begin{threeparttable}
\begin{tabular}{ ll l}
\toprule\toprule
\underline{1: Independent predictors}&\\
  $ \beta_j=\left\{\begin{array}{lll}
(-1)^{u_j}(|z_j|+4\frac{\log(n)}{\sqrt{n}}) \\
\quad j=1,\dots,5\\ 0 \text{ otherwise.}\end{array}\right.$&$\begin{array}{l}
x_{i} \sim N(0,\Sigma)\\ \Sigma=I_p\end{array}$\\
$\qquad u_j \sim \text{B}(1,0.4), z_j\sim N(0,1)$   \\ 
\midrule
\underline{2-4: Equicorrelated predictors}&\\
$\beta_{j}=\left\{\begin{array}{ll} 5 & j=1,\dots,5 \\ 0 &\text{ otherwise}
 \end{array}\right.$&
 $\begin{array}{l}
 x_{i} \sim N(0,\Sigma)\\
 \Sigma =(1-\rho)I_{p} + \rho\iota\iota'\end{array}$& $\rho=\{0.3,0.6,0.9\}$\\
\midrule
\underline{5-7: Autoregressive correlation}&\\
 $\beta_j=\left\{\begin{array}{ll}3,1.5,2 & j=1,4,7 \\ 0 &\text{ otherwise}
 \end{array}\right.$&$\begin{array}{l}
 x_{i}\sim N(0,\Sigma)\\
 \Sigma_{jk} = \rho^{|j-k|}\end{array}$ & $\rho=\{0.3,0.6,0.9\}$\\ \midrule
 \underline{8-10: Factor models}&\\
 $ \beta_j=\left\{\begin{array}{ll}5 & j=1,\dots,5 \\ 0 &\text{ otherwise}\end{array}\right.$ & $\begin{array}{l}x_{ij}=\sum_{l=1}^k f_{il}\phi_{lj}+\eta_{ij}\end{array}$ & $k=\{2,10,20\}$\\
\midrule 
\underline{11-13: Group structure}&\\
 $ \beta_j=\left\{\begin{array}{ll}3 & j=1,\dots,15 \\ 0 &\text{ otherwise}
 \end{array}\right.$&$\begin{array}{ll} 
 x_{i,k+3m}&= f_{ik}+\delta \eta_{i,k+3m}, \\
    x_{il} &=\eta_{il}\\
&\hspace{-1.5cm}k=1,2,3,\text{ and }m=0,\dots,4,\end{array} $ & $\delta^2=\{0.01,0.05,0.1\}$\\
\midrule 
\underline{14: Extreme correlation}&\\
 $ \beta_j=\begin{cases}5 & j=1,\dots,5 \\ 0 &\text{ otherwise}\end{cases} $ &\multicolumn{2}{l}{$\begin{array}{ll}
 x_{ij} & =(f_{ij}+\eta_{ij})/\sqrt{2},\\
 x_{i,j+k} &= x_{ij}+\sqrt{0.01}f_{i,j+k},\\
  x_{il}& =(f_{il}+\sum_{m=1}^5 \eta_{im})/2, \\
 &\hspace{-1cm}j = 1,\ldots,5,\text{ and } k = 5,10.\end{array}$} \\
  \bottomrule \bottomrule
    \end{tabular}
\begin{tablenotes}
\footnotesize
\item This table shows the different specifications for $\beta$ and $X$ in \eqref{eq: dgpmc} for 14 different simulation designs. Let $f_{ij}$, $\phi_{ij}$, $\eta_{ij}$ denoted standard normal random variables independent across $i$ and $j$, and $l=16,\dots,p$. All regressors are scaled by their standard deviation.
\end{tablenotes}
\end{threeparttable}
\end{table}

\subsection{Results}
\begin{figure}[tb!]
\centering
\caption{Coverage rate Monte Carlo experiments}
\includegraphics*[width=\textwidth,trim=0 0 0 0]{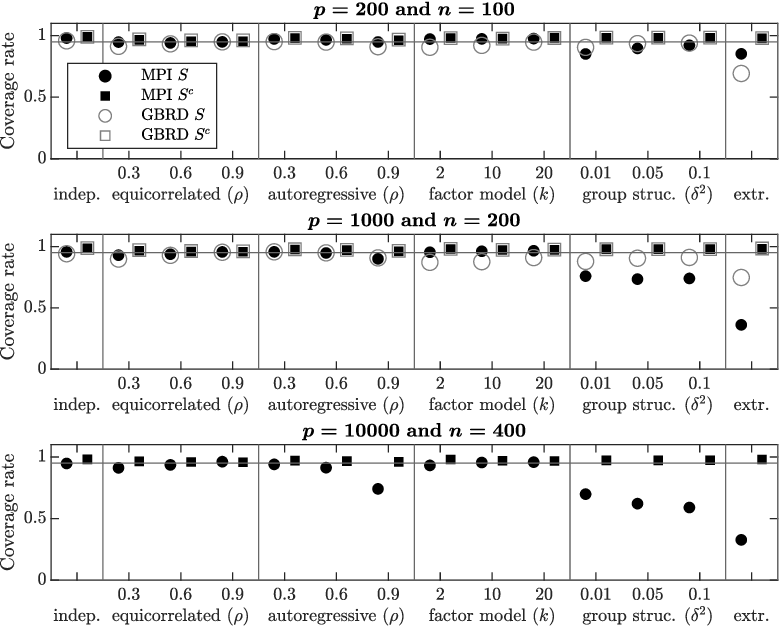}
\fnote{Note: this figure shows the coverage rates of the Moore-Penrose pseudoinverse (MPI,black markers) and the method of \citet{van2014asymptotically} (GBRD,white markers), for non-zero ($S$,circles) and zero ($S^c$,squares) coefficients. Results are based on 1000 replications of the linear model \eqref{eq: dgpmc}. The panels correspond to the dimensions $(p,n)=(200,100)$, $(p,n)=(1000,200)$, and $(p,n)=(10000,400)$, and the x-axis shows the 14 simulation designs specified in \eqref{tab:design}. 
}
\label{fig: MCcr}
\end{figure}

The results of the Monte Carlo experiments confirm the theory from \cref{sec: theory}. \cref{fig: MCcr} shows the coverage rates of MPI and GBRD for non-zero ($S$) and zero ($S^{c}$) coefficients. The coverage of MPI on $S$ and $S^{c}$ is close to nominal for all independent, equicorrelated, and factor model designs. These simulation designs satisfy the spherical symmetry assumption on the regressor matrix. For the factor model, GBRD attains a slightly lower coverage on the nonzero coefficients, and this effect increases when moving from $p=200$ to $p=1000$.

Some, but not all, of the simulation designs that violate the symmetry assumption report undercoverage for both MPI and GBRD. \cref{fig: MCcr} shows undercoverage on $S$ for autoregressive regressors with correlation parameter $\rho=0.9$, which increases with the dimensionality of the data. For both the group structure and the extreme dependence simulation designs, the undercoverage on $S$ is quite severe. 

\cref{tab:resultsmain} reports the average MAE, coverage and power for experiments with equicorrelated predictors with $\rho=0.6$, autoregressive correlation with $\rho=0.6$, and a factor model with $k=10$. The results for the remaining designs are presented in \cref{A: mcresults}. All results corroborate the undercoverage of JM on $S$ found by \citet{dezeure2015high}.

\begin{table}[t!]
  \centering \small
  \caption{Results Monte Carlo experiments with simulation design 3, 6, and 9}
  \begin{threeparttable}
    \begin{tabular}{llrrrrrrrrr}
    \toprule \toprule 
          &       & \multicolumn{3}{c}{$(p,n)=(200,100)$} & \multicolumn{3}{c}{$(p,n)=(1000,200)$} & \multicolumn{3}{c}{$(p,n)=(10000,400)$} \\
          &       & \multicolumn{1}{l}{MAE} & \multicolumn{1}{l}{CR} & \multicolumn{1}{l}{power} & \multicolumn{1}{l}{MAE} & \multicolumn{1}{l}{CR} & \multicolumn{1}{l}{power} & \multicolumn{1}{l}{MAE} & \multicolumn{1}{l}{CR} & \multicolumn{1}{l}{power} \\
         \cmidrule(lr){3-5}\cmidrule(lr){6-8}\cmidrule(lr){9-11}
method & $S$ & \multicolumn{9}{c}{Equicorrelated with $\rho=0.6$} \\
\midrule
    \multirow{2}[0]{*}{MPI} & $S$   & 3.66  & 0.94  & 0.19  & 2.11  & 0.94  & 0.43  & 1.39  & 0.94  & 0.79 \\
          & $S^c$ & 3.51  & 0.95  &       & 1.90  & 0.96  &       & 1.23  & 0.96  &  \\[+1mm]
    \multirow{2}[0]{*}{RID} & $S$   & 2.86  & 0.94  & 0.26  & 2.00  & 0.94  & 0.47  & 1.38  & 0.94  & 0.80 \\
          & $S^c$ & 2.60  & 0.96  &       & 1.76  & 0.96  &       & 1.21  & 0.96  &  \\[+1mm]
    \multirow{2}[0]{*}{GBRD} & $S$   & 2.68  & 0.93  & 0.28  & 1.96  & 0.93  & 0.48  &       &       &  \\
          & $S^c$ & 2.25  & 0.96  &       & 1.61  & 0.96  &       &       &       &  \\[+1mm]
    \multirow{2}[0]{*}{JM} & $S$   & 2.99  & 0.31  & 0.39  & 2.80  & 0.25  & 0.51  &       &       &  \\
          & $S^c$ & 0.80  & 0.95  &       & 0.61  & 0.98  &       &       &       &  \\
          \midrule
          &       & \multicolumn{9}{c}{Factor model with $k=10$} \\ \midrule
    \multirow{2}[0]{*}{MPI} & $S$   & 3.87  & 0.97  & 0.13  & 2.26  & 0.96  & 0.37  & 1.47  & 0.96  & 0.68 \\
          & $S^c$ & 3.76  & 0.98  &       & 2.10  & 0.97  &       & 1.37  & 0.97  &  \\[+1mm]
    \multirow{2}[0]{*}{RID} & $S$   & 2.57  & 0.94  & 0.27  & 1.98  & 0.94  & 0.44  & 1.44  & 0.95  & 0.70 \\
          & $S^c$ & 2.06  & 0.98  &       & 1.66  & 0.97  &       & 1.28  & 0.97  &  \\[+1mm]
    \multirow{2}[0]{*}{GBRD} & $S$   & 2.66  & 0.92  & 0.27  & 2.09  & 0.88  & 0.44  &       &       &  \\
          & $S^c$ & 1.99  & 0.98  &       & 1.45  & 0.97  &       &       &       &  \\[+1mm]
    \multirow{2}[0]{*}{JM} & $S$   & 3.40  & 0.15  & 0.43  & 3.35  & 0.10  & 0.50  &       &       &  \\
          & $S^c$ & 0.36  & 0.96  &       & 0.33  & 0.98  &       &       &       &  \\
\midrule
& & \multicolumn{9}{c}{Autoregressive with $\rho=0.6$}\\
\midrule
    \multirow{2}[0]{*}{MPI} & $S$   & 0.64  & 0.97  & 0.65  & 0.34  & 0.95  & 0.98  & 0.21  & 0.91  & 1.00 \\
          & $S^c$ & 0.61  & 0.98  &       & 0.29  & 0.97  &       & 0.17  & 0.97  &  \\[+1mm]
    \multirow{2}[0]{*}{RID} & $S$   & 0.49  & 0.95  & 0.83  & 0.31  & 0.95  & 0.99  & 0.21  & 0.91  & 1.00 \\
          & $S^c$ & 0.43  & 0.98  &       & 0.27  & 0.97  &       & 0.17  & 0.97  &  \\[+1mm]
    \multirow{2}[0]{*}{GBRD} & $S$   & 0.47  & 0.95  & 0.86  & 0.31  & 0.95  & 0.98  &       &       &  \\
          & $S^c$ & 0.39  & 0.98  &       & 0.28  & 0.97  &       &       &       &  \\[+1mm]
    \multirow{2}[0]{*}{JM} & $S$   & 0.48  & 0.67  & 0.96  & 0.33  & 0.67  & 1.00  &       &       &  \\
          & $S^c$ & 0.20  & 0.97  &       & 0.15  & 0.97  &       &       &       &  \\
          \bottomrule \bottomrule
    \end{tabular}%
    \begin{tablenotes}
\footnotesize
\item This table shows the mean absolute error of the estimated coefficients (MAE), coverage rates (CR) and statistical power of the Moore-Penrose pseudoinverse (MPI), ridge regression (RID), and the methods of \citet{van2014asymptotically} (GBRD) and \citet{javanmard2014confidence} (JM). Results are provided separately for non-zero ($S$) and zero ($S^c$) coefficients. Results are based on 1000 replications of the linear model \eqref{eq: dgpmc}, with $(p,n)=(200,100)$, $(p,n)=(1000,200)$, and $(p,n)=(10000,400)$, and coefficients and regressors specified in design 3,6, and 9 in \cref{tab:design}. For $(p,n)=(10000,400)$, only MPI and RID are included in the Monte Carlo experiment. The results for the remaining design are presented in \cref{A: mcresults}. 
\end{tablenotes}
\end{threeparttable}%
  \label{tab:resultsmain}%
  \end{table}
  
Since the equicorrelated and factor design are spherically symmetric, the coverage rates of MPI and RID are robust to the choice of $p$ and $n$. For the autoregressive model structure, we find a drop in coverage on $S$ when moving to $(p,n) = (10000,400)$. 

For all simulation designs with $p=200$ and $n=100$, the noise in the pseudoinverse translates in high MAE for MPI. This can be largely remedied by incorporating the ridge adjustment. In terms of inference, this has the benefit that the coverage rate is not distorted, while power increases. For the two higher-dimensional settings, the findings for MPI and RID are quite similar. 

\begin{figure}[tb!]
\centering
\caption{FWER of MPI and RID for equicorrelated design with $\rho=0.6$}
\includegraphics*[width=\textwidth,trim=30 95 40 95]{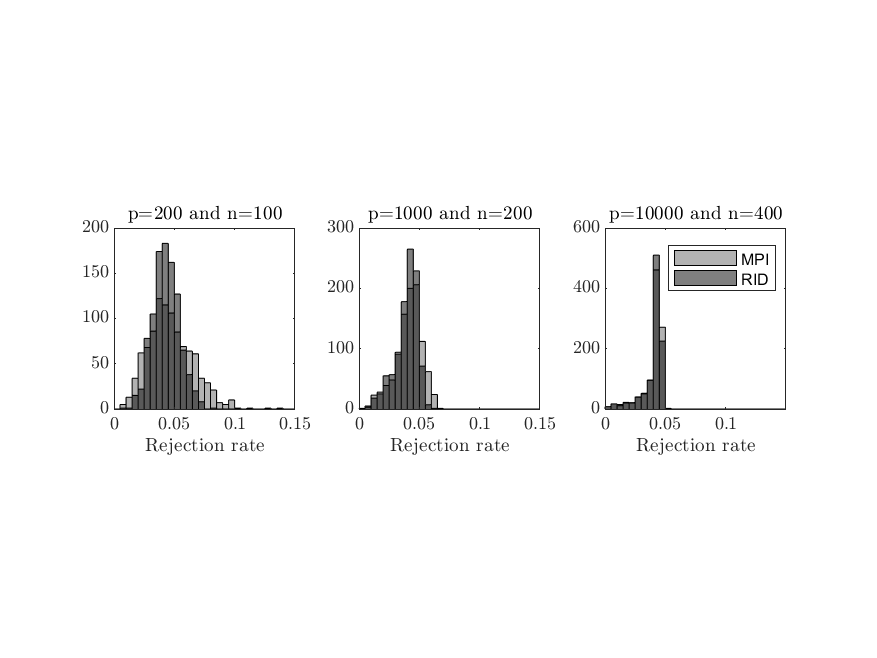}
\fnote{This figure shows the empirical distribution of the family-wise error rate (FWER) of the Moore-Penrose pseudoinverse (MPI) and ridge regression (RID) for the equicorrelated simulation design with $\rho=0.6$. For each replication, the rejection rate is calculated as the percentage of coefficient values in the data generating process that fall inside the 95\% confidence interval. The panels correspond to the dimensions $(p,n)=(200,100)$, $(p,n)=(1000,200)$, and $(p,n)=(10000,400)$. }
\label{fig: Histcr}
\end{figure}

When the number of variables is of the same order as the number of observations, the discussion below \cref{theorem:power} suggests that regularization may improve the efficiency of the confidence intervals. \cref{fig: Histcr} shows histograms of the FWER across experiments for MPI and RID under an equicorrelated design with $\rho=0.6$. In the case with $(p,n)=(200,100)$, the RID rejection rates are indeed more tightly concentrated around the nominal rate of 5\% than the MPI rejection rates. For $(p,n)=(1000,200)$ and $(10000,400)$, the number of variables is an order of magnitude larger than the number of observations and the benefit of regularization disappears. 

\begin{figure}[tb!]
\centering
\caption{Computation time inference methods for different simulation designs}
\includegraphics*[width=\textwidth,trim=30 95 40 95]{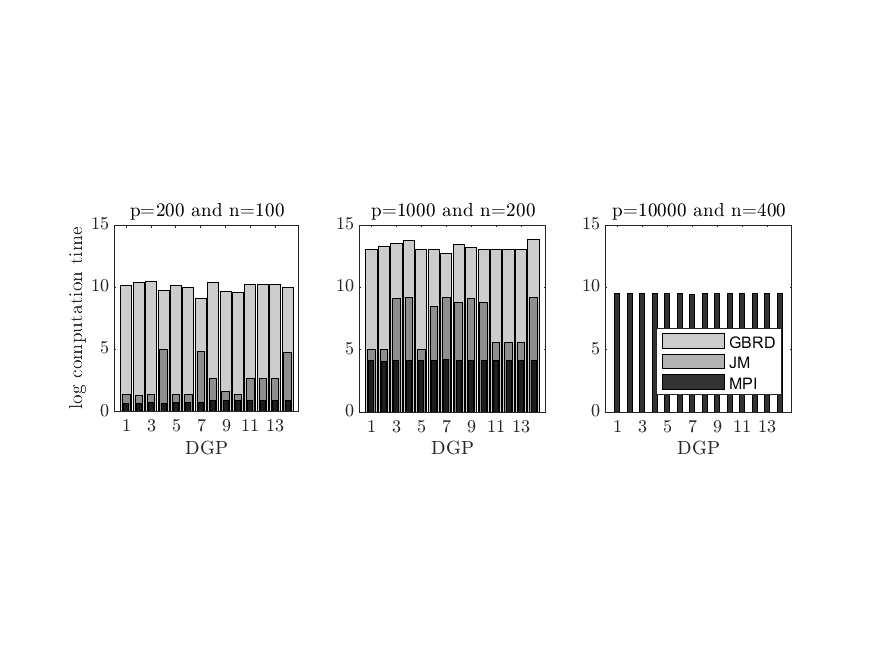}
\fnote{This figure shows the natural logarithm of the computation times of the Moore-Penrose pseudoinverse (MPI, black), and the methods of \citet{van2014asymptotically} (GBRD, light gray) and \citet{javanmard2014confidence} (JM, dark gray). The computation times are calculated as the total number of seconds required to run each inference method over 1000 replications of each simulation design. The panels correspond to the dimensions $(p,n)=(200,100)$, $(1000,200)$, and $(10000,400)$, and the x-axis corresponds to the 14 simulation designs specified in \cref{tab:design}.}
\label{fig: time}
\end{figure}

\cref{fig: time} shows the logarithm of the computation times of MPI, GBRD, and JM for each simulation design. MPI demands computation time that is negligible compared to the time needed by GBRD and JM with $p=200$ and $p=1000$. Due to the computational costs of GBRD and JM, they are not included in the experiments with $(p,n)=(10000,400)$. We find that MPI is still computationally feasible in this setting, making them scalable to ultra high-dimensional data.

\section{Conclusion}\label{sec: conclusion}
This paper proposes computationally efficient methods for constructing confidence intervals in high-dimensional linear regression models. We employ a debiasing strategy where the Moore-Penrose pseudoinverse is used as an approximate inverse of the singular empirical covariance matrix of the regressors. This strategy is shown to be valid under a symmetry assumption on the regressors. The covariance matrix of the estimates is available in closed form and free of tuning parameters. Confidence intervals can then be constructed using standard procedures. A ridge-based regularization can yield power improvements in finite samples.

We consider large scale Monte Carlo experiments that show that the proposed estimators provide valid confidence intervals with correct coverage rates, even when the symmetry assumption on the regressors does not hold.

\bibliographystyle{apalike}
\bibliography{literature}

\begin{appendices}
	\crefalias{section}{appsec}
	\crefalias{subsection}{appsubsec}
	\numberwithin{equation}{section}
		\numberwithin{nondecrl}{section}
\section{Preliminary results}\label{app:aux}
Throughout, $C$ and $C_{i}$ denote generic positive constants that can differ between lines. We abbreviate $1-C_{1}\exp(-C_{2}n)$ as $1-O(\exp(-n))$.

\subsection{Properties of spherical distributions}\label{A: propA1}
Decompose $X$ by a singular value decomposition as
\begin{equation}\label{eq:svd}
X = VSU^{\top},
\end{equation}
 where $V\in \mathcal{O}(n)$, $S$ the $n\times p$ matrix of singular values, and $U\in\mathcal{O}(p)$. Under \cref{ass:symmetry}, $X$ is invariant under right multiplication with an orthogonal matrix, and therefore $U$ is uniformly distributed on $\mathcal{O}(p)$. When $n<p$,
\begin{equation}\label{eq:svd2} 
X = VS_{n}U_{n}^{\top},
\end{equation}
where $S_{n}$ is an $n\times n$ matrix with the non-zero singular values on its diagonal, and $U_{n}$ is a $p\times n$ matrix that satisfies $
U_{n}^{\top} = [I_{n}, O_{n,p-n}]U^{\top}$.
Since $U$ is uniformly distributed over $\mathcal{O}(p)$, $U_{n}$ is uniformly distributed over the Stiefel manifold $V_{n,p}$, defined as $
V_{n,p} = \left\{A\in R^{p\times n}: A^{\top}A = I_{n}\right\}$. See for example \citet{chikuse1990matrix}.

We will use the following results.
\begin{lemma}[Lemma 5.3.2 of \citet{vershynin2018high}]\label{lem:JL} Let $v$ be a fixed $p\times 1 $ vector, and $U_{n}$ a $p\times n$ matrix that is distributed uniformly over $V_{n,p}$. Let $(C_{u},\epsilon_{u})>0$, then with probability at least  $1-2\exp(-C_{u}\epsilon_{u}^2 n)$, we have that
\begin{equation}
	(1-\epsilon_u)\frac{n}{p}\leq \frac{v^{\top}U_{n}U_{n}^{\top}v}{v^{\top}v} \leq (1+\epsilon_{u})\frac{n}{p}.
	\end{equation}
\end{lemma}

\begin{lemma}[\citet{fan2008sure}]\label{lem:fan} For any $(C_{f},t)>0$, there exists an $\epsilon_{f}>0$ such that with probability at most $3\exp(-C_{f}n) +2\frac{\exp(-t^2/2)}{t\sqrt{2\pi}}$, we have that
\begin{equation}
|e_{1}^{\top}U_{n}U_{n}^{\top}e_{2}|>(1+\epsilon_{f})\frac{t}{\sqrt{n}}\frac{n}{p}.
\end{equation}
\end{lemma}
Proof: Lemma 5 of \citet{fan2008sure} shows that for any $C_{f}>0$, there exists an $\epsilon_{f}>0$, such that
\begin{equation}\label{eq:fan}
P\left(|e_{1}^{\top}U_{n}U_{n}^{\top}e_{2}|>(1+\epsilon_{f})\frac{1}{\sqrt{n}}\frac{n}{p}|W|\right)\leq 3\exp(-C_{f}n),
\end{equation}
where $W\sim N(0,1)$ independent of $U_{n}$. For random variables $X$ and $Y$ and a constant $t$, we have $ P(X>t)\leq P(X>Y, Y>t)\leq P(X>Y)+P(Y>t)$. Set $t=(1+\epsilon_{f})\frac{1}{\sqrt{n}}\frac{n}{p}$, $X=|e_{1}^{\top}U_n U_{n}^{\top}e_{2}|/t$ and $Y=|W|$. Using \eqref{eq:fan} and the standard tail bound bound $P(|W|>t)\leq 2\frac{\exp(-t^2/2)}{t\sqrt{2\pi}}$, we have
\begin{equation} 
P\left(|e_{1}^{\top}U_{n}U_{n}^{\top}e_{2}|>(1+\epsilon_{f})\frac{t}{\sqrt{n}}\frac{n}{p}\right)\leq 3\exp(-C_{f}n) + 2\frac{\exp(-t^2/2)}{t\sqrt{2\pi}}.
\end{equation}

\begin{lemma}\label{lem:use} Let $\hat{\Sigma}_{ii} = e_{i}^{\top}\frac{1}{n}X^{\top}Xe_{i}$. Under \cref{ass:factor} there exists constants $0<c<C$, such that $c<\hat{\Sigma}_{ii}<C$ with probability  $1-O(\exp(-n))$.
\end{lemma}
Proof: Notice that
\begin{equation}
\begin{split}
    P(\hat{\Sigma}_{ii}<c)&\leq P\left(\frac{1}{n}e_{i}^{\top}R^{\top}Re_{i}-\frac{2}{n}|e_{i}^{\top}R^{\top}F\Lambda e_{i}|<c\right)\\
    &\leq \underbrace{P\left(\frac{1}{n}e_{i}^{\top}R^{\top}Re_{i}<c+\tilde{c}\right)}_{(I)} + \underbrace{P\left(\frac{2}{n}|e_{i}^{\top}R^{\top}F\Lambda^{\top} e_{i}|>\tilde{c}\right)}_{(II)}.
    \end{split}
\end{equation}
We first bound (I). By \cref{ass:R}, $Re_{i}$ is a subgaussian vector with independent elements. Under this assumption, $u_{j}=e_{i}^{\top}R^{\top}e_{j}e_{j}^{\top}Re_{i}-\Sigma_{ii}^{R}$ is subexponential with mean zero. From \cref{ass:R}, we also have that $\Sigma_{ii}^{R}=e_{i}^{\top}\Sigma^{R}e_{i}\geq \lambda_{\min}(\Sigma^{R})>0$. Pick $c,\tilde{c}>0$, such that $\lambda_{\min}(\Sigma^R)-(c+\tilde{c})=c_{1}>0$. Then, 
\begin{equation}
    \begin{split}
        (I) &\leq  P\left(\frac{1}{n}\sum_{j=1}^{n}u_{j} \geq c_{1}\right)=O(-\exp(-n)).
    \end{split}
\end{equation}
Continuing with (II), we first condition on $\Lambda$. From \cref{ass:factor}, $e_{j}F\Lambda^{\top}e_{i}$ is a subgaussian random variable, and by independence between $R$ and $(F,\Lambda)$, $v_{jl}=e_{i}^{\top}R^{\top}e_{j}e_{j}^{\top}Fe_{l}e_{l}^{\top}\Lambda^{\top}e_{i}$ where $j=1,\ldots,n$ and $l=1,\ldots,k$, is subexponential with mean zero. Then,
\begin{equation}
    \begin{split}
        P\left(\left.\frac{2}{n}|\sum_{j=1}^{n}e_{i}^{\top}R^{\top}e_{j}e_{j}^{\top}Fe_{l}e_{l}^{\top}\Lambda^{\top} e_{i}|>\tilde{c}\right|\Lambda\right)&=O\left(\exp\left(-\frac{n}{|e_{l}^{\top}\Lambda^{\top} e_{i}|}\right)\right).
    \end{split}
\end{equation}
 Taking the expectation over $\Lambda$, using Jensen's inequality and the fact that $\text{E}[|e_{l}^{\top}\Lambda^{\top} e_{i}|]\leq C<\infty$,
\begin{equation}
    \begin{split}
        P\left(\frac{2}{n}|e_{i}^{\top}R^{\top}Fe_{l}e_{l}^{\top}\Lambda^{\top} e_{i}|>\tilde{c}\right)&=O(\exp(-n)).
    \end{split}
\end{equation}
Since $l=1,\ldots,k$ with $k$ finite, we then also have that 
\begin{equation}\label{eq:needed}
    \begin{split}
        P\left(\frac{2}{n}|e_{i}^{\top}R^{\top}F\Lambda^{\top} e_{i}|>\tilde{c}\right)&=O(\exp(-n)).
    \end{split}
\end{equation}
What remains to be shown is that there exists a $C$ such that 
\begin{equation}
    P\left(\hat{\Sigma}_{ii}>C\right)=O(\exp(-n)).
\end{equation}
For this it is sufficient to show that 
\begin{equation}
    \begin{split}
    P\left(\frac{1}{n}\Lambda_{i}^{\top} F^{\top}F\Lambda_{i}>C_1\right)= O(\exp(-n)),\quad 
    P\left(\frac{1}{n}e_{i}^{\top} R^{\top}Re_{i}>C_2\right)=O(\exp(-n)), \notag
\end{split}
\end{equation}
where $\Lambda$ has rows $\Lambda_{i}$. Starting with the first, notice that
\begin{equation}
    \begin{split}
    P\left(\frac{1}{n}\Lambda_{i}^{\top} F^{\top}F\Lambda_{i}>C_1\right)
    &=P\left(\sum_{l,m=1}^{k}\Lambda_{il}\Lambda_{im}\frac{1}{n}f_{l}^{\top}f_{m}>C_1\right)\\
    &\leq P\left(\frac{1}{2}\sum_{l,m=1}^{k}\Lambda_{il}^2\frac{1}{n}f_{l}^{\top}f_{l} + \Lambda_{im}^2\frac{1}{n}f_{m}^{\top}f_{m}>C_1\right)\\
    &\leq P\left(k^2\max_{l=1,\ldots,k}\Lambda_{il}^2\frac{1}{n}f_{l}^{\top}f_{l}>C_1\right)\\
    &\leq k P\left(k^2\Lambda_{il}^2\frac{1}{n}f_{l}^{\top}f_{l}>C_1\right)
\end{split}
\end{equation}
Conditioning on $\Lambda_{il}$ and using that the elements of $f_{l}$ are independent subgaussian,
\begin{equation}
    \begin{split}
        P\left(\Lambda_{il}^2\frac{1}{n}f_{l}^{\top}f_{l}>C_1|\Lambda_{il}\right)\leq C_{3}\exp\left(-C_{4}\frac{n}{\Lambda_{il}^2}\right).
    \end{split}
\end{equation}
Taking the expectation over $\Lambda_{il}$, using Jensen's inequality and the fact that $\text{E}[\Lambda_{il}^2]\leq C<\infty$, we then find that for fixed $k$, there exists a constant $C_1>0$ such that
\begin{equation}
    P\left(\frac{1}{n}\Lambda_{i}^{\top} F^{\top}F\Lambda_{i}>C_1\right)=O(\exp(-n)).
\end{equation}
By the same arguments as in Section 6.2.2 of \citet{javanmard2014confidence}, we also have that there exists a constant $C_2>0$, such that
\begin{equation}
    P\left(\frac{1}{n}e_{i}^{\top}R^{\top}Re_{i}>C_2\right)=O(\exp(-n)).
\end{equation}
This completes the proof.

\begin{lemma}\label{lem:mineig} Suppose \cref{ass:sparsity} and \cref{ass:factor} hold. Denote by $\hat{\lambda}_{\min\neq 0}$ the smallest nonzero eigenvalue of $p^{-1}X^{\top}X$. Then for some constant $C>0$, $\hat{\lambda}_{\min\neq 0}\geq C$ with probability $1-O(\exp(-n))$.
\end{lemma}

To determine the minimum nonzero eigenvalue, define $v$ such that $v^{\top}v=1$ and write
\begin{equation}\label{eq:C25}
\begin{split}
    v^{\top}p^{-1}X^{\top}Xv& = v^{\top}p^{-1}\Lambda F^{\top}F\Lambda^{\top}v + v^{\top}p^{-1}R^{\top}Rv + 2p^{-1}v^{\top}\Lambda F^{\top}Rv \\
    &\geq v^{\top}p^{-1}\Lambda F^{\top}F\Lambda^{\top}v + v^{\top}p^{-1}R^{\top}Rv -2p^{-1}|v^{\top}\Lambda F^{\top}Rv| \\
    &\geq v^{\top}p^{-1}\Lambda F^{\top}F\Lambda^{\top}v + v^{\top}p^{-1}R^{\top}Rv - C_{1}(n/p)\\
    &\geq v^{\top}p^{-1}R^{\top}Rv - C_{1}(n/p).
\end{split}
\end{equation}
The second to last line holds uniformly over $v$ with probability $1-O(\exp(-n))$ by \eqref{eq:needed} replacing $e_{i}$ by $v$ and noting in applying Jensen's inequality that $\Lambda$ has subgaussian columns by \cref{ass:factor}. This shows that $\lambda_{\min\neq 0}(p^{-1}X^{\top}X)\geq \lambda_{\min\neq 0}(p^{-1}R^{\top}R)-C_{1}(n/p)$. We will now show that the minimum nonzero eigenvalue of $p^{-1}R^{\top}R$ is bounded away from zero by a constant. This implies that for $p$ sufficiently larger than $n$, \cref{lem:mineig} holds. 

Since $p>n$, the minimum nonzero eigenvalue of $p^{-1}R^{\top}R$ is equal to  \begin{equation}
\begin{split}
    \lambda_{\min\neq 0}(p^{-1}R^{\top}R)&=\lambda_{\min}(p^{-1}R R^{\top}) \\
    &= \lambda_{\min}(p^{-1}\tilde{R}\Sigma^{R}\tilde{R})\\
    &\geq \lambda_{\min}(\Sigma^{R})\lambda_{\min}(p^{-1}\tilde{R}\tilde{R}^{\top})\\
    &\geq C\lambda_{\min}(p^{-1}\tilde{R}\tilde{R}^{\top}),
    \end{split}
    \end{equation}
    with $\tilde{r}_{i}=e_{i}^{\top}\tilde{R}$ defined in \cref{ass:R} and the last line follows from \cref{ass:R}. Since $\tilde{R}$ is an $n\times p$ matrix with independent subgaussian elements, by \citet{vershynin2010introduction} Theorem 5.39,  $\lambda_{\min}(p^{-1}\tilde{R}\tilde{R}^{\top})\geq 1-C_{2}\sqrt{n/p}$ with probability $1-O(\exp(-n))$. 
    
    We now have shown that $\lambda_{\min\neq 0}(p^{-1}X^{\top}X)\geq 1-C_{2}\sqrt{n/p}-C_{1}(n/p)$. From \cref{ass:sparsity}, we can assume that $p>\max(nC_{2}^{2}\epsilon^{-2},nC_{1}\epsilon^{-1})$ for some positive constant $\epsilon<1$. Then with probability $1-O(\exp(-n))$, $\hat{\lambda}_{\min\neq 0}(p^{-1}X^{\top}X)\geq C>0$. 

\section{Proofs of the main theory}

\subsection{Proof of \cref{theorem:main}}\label{app:prooftheoremmain}
Rewrite the estimator $\hat{\beta}$ in \eqref{eq:betac} as
\begin{equation}\label{eq: corr}
\begin{split}
\sqrt{n}\left(\hat{\beta}-\beta\right) &= \Delta+Z\\
\Delta &= \sqrt{n}\left(\hat{M}X-I_{p}\right)\left(\beta-\hat{\beta}^{\text{lasso}}\right)\\
Z&=\sqrt{n}\hat{M}\varepsilon.
\end{split}
\end{equation}

The bias term of the estimator in \eqref{eq: corr} is bounded by the norm inequality,
\begin{equation}\label{eq:splitbias}
||\Delta||_{\infty}\leq \sqrt{n}\left|\left|\hat{M}X-I_{p}\right|\right|_{\max}||\beta-\hat{\beta}^{\text{lasso}}||_{1}.
\end{equation}
We then have that
\begin{equation}
    \begin{split}
        P\left(||\Delta||_{\infty}\geq a\frac{s_{0}\log p}{\sqrt{n}}\right)&\leq P\left(||\hat{M}X-I_{p}||_{\max}\geq \sqrt{\frac{\log p}{n}}b^{-1}a\right)\\
        &\quad + P\left(||\beta-\hat{\beta}^{\text{lasso}}||_{1}\geq bs_{0}\sqrt{\frac{\log p}{n}}\right).
    \end{split}
\end{equation}
\cref{theorem:biasvanish} provides a probability bound for $\left|\left|\hat{M}X-I_{p}\right|\right|_{\max}$, and \cref{lem:acclasso} provides a probability bound for $||\beta-\hat{\beta}^{\text{lasso}}||_{1}$.

\begin{lemma}\label{theorem:biasvanish} \sloppy Suppose \cref{ass:sparsity} and \ref{ass:symmetry} holds. Define $\hat{M}$ as in \eqref{eq: mppi} and $\hat{D}$ as in \eqref{eq: diagmppi}. Then, there exists an $c>0$ such that
\begin{equation}\label{eq:Pbound}
P\left(\left|\left|\hat{M}X-I_{p}\right|\right|_{\max}\geq c\sqrt{\frac{\log p}{n}}\right) = O(p^{-1}).
\end{equation}
\end{lemma}
Proof: \cref{app:prooftheorembiasvanish}. 

\begin{lemma}\label{lem:acclasso}
 Suppose \cref{ass:sparsity}, \ref{ass:error} and either \ref{ass:R} or \ref{ass:factor} hold. Consider the lasso estimator \eqref{eq:lasso} with $\lambda\geq 4C_{1}\sigma\sqrt{\log p/n}$ for some sufficiently large constant $C_{1}$,
\begin{equation}\label{eq:accguar}
\begin{split}
&P\left(\left|\left|\beta-\hat{\beta}^{\text{lasso}}\right|\right|_{1}\geq b s_{0}\sqrt{\frac{\log p}{n}}\right)=O(p^{-1}).
\end{split}
\end{equation}
\end{lemma}
Proof: \cref{proof:compcond}. 

Combining \cref{ass:sparsity}, \cref{theorem:biasvanish}, and \cref{lem:acclasso}, it follows that the bias can be bounded by
\begin{equation}\label{eq:biasvanish2}
    ||\Delta||_{\infty} = O_{p}\left(s_{0}\frac{\log p}{\sqrt{n}}\right)= o_{p}(1).    
\end{equation}

From \cref{ass:error} and \eqref{eq: corr} it follows that $Z|X\sim N(0,\hat{\Omega})$ with $\hat{\Omega} = \hat{\sigma}^2n\hat{M}\hat{M}^\top$. The following lemma shows that $\hat{\Omega}_{jj} = O_{p}(1)$, which completes the proof of \cref{theorem:main}.

\begin{lemma}\label{theorem:finitevar} Suppose \cref{ass:sparsity} -- \ref{ass:error} hold. Define $ Z_{j} = \sqrt{n}\hat{D}_{jj}x_{j}^\top (XX^\top )^{-1}\varepsilon$ and $\hat{\Omega}_{jj} = n\hat{\sigma}^2 e_{j}^\top \hat{M}\hat{M}^\top e_{j}$ with  $\hat{M}$ as before. For $j=1,\ldots,p$,
\begin{equation}
    Z_{j}|X \sim N(0,\hat{\Omega}_{jj}),\quad
    \hat{\Omega}_{jj}= O_{p}(1).
\end{equation}
\end{lemma}
Proof: \cref{app:prooftheoremfinitevar}. 

\subsection{Proof of \cref{theorem:main2}}\label{app:prooftheoremmain2}
We show that for a sufficiently small penalty parameter $\lambda$, the results under a Moore-Penrose inverse carry over to a ridge adjusted estimator. We adjust \cref{theorem:biasvanish} and \cref{theorem:finitevar} to apply to the ridge adjusted estimator. \cref{theorem:main2} then follows.

\begin{lemma}\label{theorem:biasvanishridge} \sloppy Suppose \cref{ass:sparsity} and \ref{ass:symmetry} holds. Define $\hat{M}$ and $\hat{D}$ as in \eqref{eq: ridge} and \eqref{eq: diagridge}. Choose $\gamma = O\left(p\sqrt{\frac{\log p}{n}}\right)$. Then, there exists an $a>0$ such that
\begin{equation}\label{eq:Pboundridge}
P\left(\left|\left|\hat{M}X-I_{p}\right|\right|_{\max}\geq a\sqrt{\frac{\log p}{n}}\right) = O(p^{-1}).
\end{equation}
\end{lemma}
Proof: \cref{app:proofridgelemma}.

\begin{lemma}\label{theorem:finitevarridge} Suppose \cref{ass:sparsity} -- \ref{ass:error} hold. Define $ Z_{j} = \sqrt{n}\hat{D}_{jj}e_{j}^\top (X^\top X + \gamma I_{p})^{-1}X^\top \varepsilon$ and $\hat{\Omega}_{jj} = n\hat{\sigma}^2 e_{j}^\top \hat{M}\hat{M}^\top e_{j}$ with  $\hat{M}$ as before. For $j=1,\ldots,p$,
\begin{equation}
    Z_{j}|X \sim N(0,\hat{\Omega}_{jj}),\quad
    \hat{\Omega}_{jj}= O_{p}(1).
\end{equation}
\end{lemma}
Proof: \cref{app:prooftheoremfinitevarridge}. 

\subsection{Proof of \cref{theorem:power}}\label{app:prooftheorempower}
The relative efficiency is determined by the difference 
\begin{equation}
    \frac{\hat{\Omega}_{ii}^{\text{MPI}}(\hat{D})-\hat{\Omega}_{ii}^{\text{RID}}(\hat{D})}{n^2 \hat{D}_{ii}^2}.
\end{equation} 
   with $\hat{\Omega}_{ii}^{\text{MPI}}$ as in \cref{theorem:main} and $\hat{\Omega}_{ii}^{\text{RID}}$ as in \cref{theorem:main2}, but both estimators are assumed to either use $\hat{D}$ as in \eqref{eq: diagmppi} or $\hat{D}$ as in \eqref{eq: diagridge}. Define the $p\times p$ matrix ${Q}$ as a diagonal matrix with containing the eigenvalues of $X^{T}X$ and by ${Q}_{+}^{-1}$ the Moore-Penrose inverse of ${Q}$. For the ridge regularized inverse we have
\begin{equation}
\begin{split}
  \hat{\Omega}_{ii}^{\text{RID}}(\hat{D})/(n^2 \hat{D}_{ii}^2) 
  & = e_{i}^{\top}(X^{\top}X+ \gamma I_{p})^{-1}X^{\top}X(X^{\top}X+\gamma I_{p})^{-1}e_{i}\\
    & = e_{i}^{\top}{U}_{n}({Q} + \gamma I_{p})^{-2}{Q}{U}_{n}^{\top}e_{i}\\
    &= e_{i}^{\top}{U}_{n}{Q}^{-1/2}_{+}A_{\text{RID}}^{2}{Q}^{-1/2}_{+}{U}_{n}^{\top}e_{i},
\end{split}
\end{equation}
where $A_{\text{RID}}^2$ is a diagonal matrix with the diagonal elements satisfying $0< A_{ii}^2< 1$. 

For the Moore-Penrose pseudoinverse,
\begin{equation}
\begin{split}
   \hat{\Omega}_{ii}^{\text{MPI}}(\hat{D})/(n^2 \hat{D}_{ii}^2)& = e_{i}^{\top}X^{\top}(XX^{\top})^{-2}Xe_{i}\\
    & = e_{i}^{\top}{U}_{n}{Q}^{-1}_{+}{U}_{n}^{\top}e_{i}\\
    &= e_{i}^{\top}{U}_{n}{Q}^{-1/2}_{+}A_{\text{MPI}}^{2}{Q}^{-1/2}_{+}{U}_{n}^{\top}e_{i},
\end{split}
\end{equation}
where $A_{\text{MPI}}^2=I$ is a diagonal matrix with the diagonal elements $A_{ii}^2=1$. 

\section{Proofs of the lemmas}\label{app:proofs}

\subsection{Proof of \cref{theorem:biasvanish}}\label{app:prooftheorembiasvanish}
As the diagonal elements of $\hat{M}X-I_{p}$ are identically zero, we need to establish a probability bound on the off-diagonal elements of $\hat{M}X-I_{p}$. These elements are of the form $e_{i}^{\top}U_{n}U_{n}^{\top}e_{j}/e_{i}^{\top}U_{n}U_{n}^{\top}e_{i}$. Fix $i\neq j$. Then,
\begin{equation}
\begin{split}
    &P\left(\frac{|e_{i}^{\top}U_{n}U_{n}^{\top}e_{j}|}{e_{i}^{\top}U_{n}U_{n}^{\top}e_{i}}>\frac{1+\epsilon_{f}}{1-\epsilon_{u}}\frac{t}{\sqrt{n}}\right)\\
    &\leq  P\left(|e_{i}^{\top}U_{n}U_{n}^{\top}e_{j}|>(1+\epsilon_{f})\frac{t}{\sqrt{n}}\frac{n}{p}\right)+ P\left(e_{i}^{\top}U_{n}U_{n}^{\top}e_{i}<(1-\epsilon_{u})\frac{n}{p}\right)\\
    &\leq 3\exp(-C_{f}n)+2\exp(-C_{u}\epsilon_{u}^2 n)+2\frac{\exp(-t^2/2)}{t\sqrt{2\pi}},
    \end{split}
\end{equation}
where we use \cref{lem:JL} and \cref{lem:fan}.

Taking the union bound over all $i\neq j$, we get
\begin{equation}
\begin{split}
    &P\left(\max_{i\neq j}\frac{|e_{i}^{\top}U_{n}U_{n}^{\top}e_{j}|}{e_{i}^{\top}U_{n}U_{n}^{\top}e_{i}}>\frac{1+\epsilon_{f}}{1-\epsilon_{u}}\frac{t}{\sqrt{n}}\right)\\
    &\leq 3\exp(-C_{f}n+2\log p)+2\exp(-C_{u}\epsilon_{u}^2n+2\log p)+2\frac{\exp(-t^2/2+2\log p)}{t\sqrt{2\pi}}.\notag
    \end{split}
\end{equation}
Now choose $t=\frac{1-\epsilon_{u}}{1+\epsilon_{f}}c\sqrt{\log p}$. Suppose $c>3\frac{1+\epsilon_{f}}{1-\epsilon_{u}}$. According to \cref{ass:sparsity}, $n>3\max(C_{f}^{-1},C_{u}^{-1}\epsilon_{u}^{-2})\log p$. Then,
\begin{equation}
\begin{split}
    P\left(\max_{i\neq j}\frac{|e_{i}^{\top}U_{n}U_{n}^{\top}e_{j}|}{e_{i}^{\top}U_{n}U_{n}^{\top}e_{i}}>c\sqrt{\frac{\log p}{n}}\right)&=O(p^{-1}).
    \end{split}
\end{equation}

\subsection{Proof of \cref{lem:acclasso}}\label{proof:compcond}
For the bound in \cref{lem:acclasso} we first show that $\max_{i=1,\ldots,p}|e_{i}^{\top}X^\top\varepsilon|/\sqrt{n\sigma^2}\leq C_{1}\sqrt{\log p}$ with high probability. Denoting by $S_{0}$ the index set of nonzero coefficients, this means that with high probability the lasso estimator satisfies $||\hat{\beta}_{S_{0}^{c}}||_{1}\leq 3||\hat{\beta}_{S_{0}}-\beta_{S_0}||_1$ by Lemma 6.3 of \citet{buhlmann2011statistics}. 

We start by decomposing
\begin{equation}
    \begin{split}
        \max_{i=1,\ldots,p}|e_{i}^{\top}X^{\top}\varepsilon/\sqrt{n\sigma^2}| \leq  \underbrace{\max_{i=1,\ldots,p}|e_{i}^{\top}\Lambda F^{\top} \varepsilon/\sqrt{n\sigma^2}|}_{(I)} +\underbrace{ \max_{i=1,\ldots,p}|e_{i}^{\top}R^{\top}\varepsilon/\sqrt{n\sigma^2}|}_{(II)}. 
    \end{split}
\end{equation}
First consider (II). Denote by $Z_{i}$ a standard normal random variable and $\hat{\Sigma}_{ii}=e_{i}^{\top}\frac{1}{n}R^{\top}Re_{i}$. By \cref{ass:error}, $\varepsilon\sim N(0,\sigma^2 I_{n})$, and (II) is distributed as
\begin{equation}
\max_{i=1,\ldots,p}\sqrt{\hat{\Sigma}_{ii}}|Z_{i}|\leq \sqrt{\max_{i=1,\ldots,p}\hat{\Sigma}_{ii}}\max_{i=1,\ldots,p}|Z_{i}|.
\end{equation}

This term can be bounded in probability as
\begin{equation}
\begin{split}
    P\left(\max_{i=1,\ldots,p}\sqrt{\hat{\Sigma}_{ii}}|Z_{i}|\geq a\sqrt{\log p}\right) &=p\left[P\left(\left.|Z_{i}|\geq a \sqrt{\hat{\Sigma}_{ii}^{-1}\log p}\right|\hat{\Sigma}_{ii}\leq C\right)P(\hat{\Sigma}_{ii}\leq C) \right.\\
    &\quad\left.+P\left(\left.|Z_{i}|\geq a \sqrt{\hat{\Sigma}_{ii}^{-1}\log p}\right|\hat{\Sigma}_{ii}> C\right)P(\hat{\Sigma}_{ii}> C)\right]\\
    &= p \left[P\left(|Z_{i}|\geq a\sqrt{C^{-1}\log p}\right) +O(\exp(-n))\right]\\
 &\leq  C_{1}\exp\left(-\frac{a^2\log p}{2C}+\log p\right)+O(p\exp(-n)). \notag
 \end{split}   
\end{equation}
where we first use \cref{lem:use}, and subsequently the tail bound for a standard normal random variable.

Choosing $a = 2\sqrt{C}$, the first term on the left hand-side is $O(p^{-1})$, while the second term is $O(\exp(-n))$ when $n$ is sufficiently larger than $\log p$ (\cref{ass:sparsity}). We conclude that,
\begin{equation}
    P\left(\max_{i=1,\ldots,p}|e_{i}^{\top}R^{\top}\varepsilon|/\sqrt{n\sigma^2}\geq a\sqrt{\log p}\right)= O(p^{-1}).
\end{equation}

For (I), notice that 
\begin{equation}
    \max_{i=1,\ldots,p}\left|\frac{1}{\sqrt{n\sigma^2}}e_{i}^{\top}\Lambda F^{\top}\varepsilon\right|\leq\max_{i=1,\ldots,p}k\max_{j=1,\ldots,k}\left|\Lambda_{ij}\right|\frac{1}{\sqrt{n\sigma^2}}\left|f_{j}^{\top}\varepsilon\right|.
\end{equation}

Define $\hat{\Sigma}^{F}_{jj} = e_{j}^{\top}\frac{1}{n}F^{\top}Fe_{j}$ and note again that $|f_{j}^{\top}\varepsilon|/\sqrt{n\sigma^2}$ is distributed as $\sqrt{\hat{\Sigma}^{F}_{jj}}|Z_{j}|$. Conditional on  $\Lambda$, we can use the same arguments as before to show that 
\begin{equation}
\begin{split}
    &P\left(\left.\max_{i=1,\ldots,p}\max_{j=1,\ldots,k}k|\Lambda_{ij}|\sqrt{\hat{\Sigma}^{F}_{jj}}|Z_{j}|\geq a\sqrt{\log p}\right|\Lambda\right)\\
 &\leq  k\exp\left(-\frac{a^2\log p}{2k^2\Lambda_{ij}^2(1+\epsilon)^2}+\log p\right)+k\exp(-Cn+\log p).
 \end{split}   
\end{equation}
Taking now the expectation of $\Lambda$ and using Jensen's inequality, we obtain
\begin{equation}
\begin{split}
    &P\left(\max_{i=1,\ldots,p}\max_{j=1,\ldots,k}k|\Lambda_{ij}|\sqrt{\hat{\Sigma}^{F}_{jj}}|Z_{j}|\geq a\sqrt{\log p}\right)\\
 &\leq  k\exp\left(-\frac{a^2\log p}{2k^2\text{E}[\Lambda_{ij}^2](1+\epsilon)^2}+\log p\right)+k\exp(-Cn+\log p).
 \end{split}   
\end{equation}
We conclude that for sufficiently large $a$ and $n$, we have that
\begin{equation}
    P\left(\max_{i=1,\ldots,p}|e_{i}^{\top}\Lambda F^{\top}\epsilon|/\sqrt{n}\sigma^2\geq a\sqrt{\log p}\right)\leq O(p^{-1}).
\end{equation}

From these results, we conclude that there exists an $C_{1}$, such that with high probability
\begin{equation}
    \max_{i=1,\ldots,p}|e_{i}^{\top}X^{\top}\varepsilon|/\sqrt{n\sigma^2}\leq C_{1}\sqrt{\log p}.
\end{equation}
Choosing then the lasso penalty parameter as $\lambda = 2C_{1}\sigma \sqrt{\frac{\log p}{n}}$ for some sufficiently large constant $C_{1}$, the results in \citet{buhlmann2011statistics} show that with high probability $||\hat{\beta}_{S_{0}^{c}}||_{1}\leq 3||\hat{\beta}_{S_{0}}-\beta_{S_{0}}||_{1}$. 

The standard $\ell_1$ norm accuracy bounds for the lasso follow if the following restricted eigenvalue condition holds \citep{bickel2009simultaneous}.
\begin{definition}\label{def:compatibility}
	The restricted eigenvalue condition holds if
	\begin{equation}\label{eq: comp}
	||\beta_{S_{0}}||_{2}\leq \frac{\||X\beta||_{2}}{\sqrt{n}\phi_{0}},
	\end{equation}
	for all $\beta$ for which $||\beta_{S_{0}^{c}}||_{1}\leq 3||\beta_{S_{0}}||_{1}$ and $\phi_{0}>0$. 
\end{definition}

We start with the following lower bound.
\begin{equation}\label{eq:a9}
\begin{split}
\frac{\beta^{\top}\frac{1}{n}X^{\top}X\beta}{\beta^{\top}\beta}
&\geq \frac{\beta^{\top}\frac{1}{n}R^{\top}R\beta}{\beta^{\top}\beta} + 2\frac{\beta^{\top}\Lambda\frac{1}{n}F^{\top}R\beta}{\beta^{\top}\beta}\\
&\geq \underbrace{\frac{\beta^{\top}\frac{1}{n}R^{\top}R\beta}{\beta^{\top}\beta}}_{(I)} - \underbrace{2\sqrt{\lambda_{\max}(\Sigma^{R})}\max_{j=1,\ldots,k}\left|\frac{\beta^{\top}\Lambda_{j}}{\sqrt{\beta^{\top}\beta}}\right|\max_{j=1,\ldots,k}\left|\frac{\frac{1}{n}f^{\top}_{j}\tilde{R}\tilde{\beta}}{\sqrt{\tilde{\beta}^{\top}\tilde{\beta}}}\right|}_{(II)},
\end{split}
\end{equation}
where $\tilde{R}$ has rows $\tilde{r}_{i}$ defined in \cref{ass:R}, $\tilde{\beta}=\Sigma^{1/2}\beta$, and $\Lambda_{j}$ denotes the $j$th column of $\Lambda$. 
(I) satisfies \cref{def:compatibility} with probability exceeding $1-O(\exp(-n))$ under \cref{ass:R} by the proof in Section 6.2.1 of \citet{javanmard2014confidence}. For (II), since $\Lambda$ and $F$ have identically distributed columns,
\begin{equation}\label{eq:c15}
    \begin{split}
        P\left(\left.\max_{j=1,\ldots,k}\left|\frac{\beta^{\top}\Lambda_{j}}{\sqrt{\beta^{\top}\beta}}\right|\max_{j=1,\ldots,k}\left|\frac{\frac{1}{n}f^{\top}_{j}\tilde{R}\tilde{\beta}}{\sqrt{\tilde{\beta}^{\top}\tilde{\beta}}}\right|>a\right|\Lambda\right)\leq k^2 P\left(\left.\left|\frac{\beta^{\top}\Lambda_{j}}{\sqrt{\beta^{\top}\beta}}\right|\left|\frac{\frac{1}{n}f^{\top}_{j}\tilde{R}\tilde{\beta}}{\sqrt{\tilde{\beta}^{\top}\tilde{\beta}}}\right|>a\right|\Lambda\right)\\
        \leq C_{3}\exp\left(-C_{4}a\frac{\sqrt{\beta^{\top}\beta}}{|\beta^{\top}\Lambda_{j}|}n\right)\notag
    \end{split}
\end{equation}
Taking the expectation over $\Lambda_{j}$, using Jensen's inequality and the fact that $\Lambda$ has subgaussian columns by \cref{ass:factor}, we have that, uniformly over $\beta$,
\begin{equation}
     P\left(\max_{j=1,\ldots,k}\left|\frac{\beta^{\top}\Lambda_{j}}{\sqrt{\beta^{\top}\beta}}\right|\max_{j=1,\ldots,k}\left|\frac{\frac{1}{n}f^{\top}_{j}\tilde{R}\tilde{\beta}}{\sqrt{\tilde{\beta}^{\top}\tilde{\beta}}}\right|>a\right)=O(\exp(-an)).
\end{equation}
Since under \cref{ass:factor}, $\lambda_{\max}(\Sigma^{R})$ is finite, we conclude that $(II)$ is $O(\log p/n)$ with probability $1-O(p^{-1})$. 

\subsection{Proof of \cref{theorem:finitevar}}\label{app:prooftheoremfinitevar}
Normality follows from \cref{ass:error}. What remains is to show that $\hat{\Omega}_{jj}=O_{p}(1)$.
\begin{equation}\label{eq:varOmega}
    \hat{\Omega}_{jj} =\text{var}(Z_{i}|X)=\sigma^2 \frac{\frac{1}{n}x_{j}^{\top}(\frac{1}{p}XX^{\top})^{-2}x_{j}}{(\frac{1}{n}x_{j}^{\top}(\frac{1}{p}XX^{\top})^{-1}x_{j})^2}.
\end{equation}
By Cauchy-Schwarz, 
\begin{equation}
    \left(\frac{1}{n}x_{j}^{\top}\left(\frac{1}{p}XX^{\top}\right)^{-1}x_{j}\right)^2\leq \left(\frac{1}{n}x_{j}^{\top}\left(\frac{1}{p}XX^{\top}\right)^{-2}x_{j}\right)\left(\frac{1}{n}x_{i}^{\top}x_{j}\right).
\end{equation}
It follows that
\begin{equation}
    \hat{\Omega}_{jj}\geq \sigma^2 (n^{-1}x_{j}^{\top}x_{j})^{-1}=\sigma^2(e_{j}^{\top}n^{-1}X^{\top}Xe_{j})^{-1}.
\end{equation}
By \cref{lem:use}, with probability at least $1-O(\exp(-n))$, we have $\hat{\Omega}_{jj}\geq C^{-1}\sigma^2$. Taking a union bound over $j=1,\ldots,p$, we find that for $n$ sufficiently larger than $\log p$, $\min_{i=1,\ldots,p}\hat{\Omega}_{jj}\geq C_{1}$ for some positive constant $C_{1}$ with probability $1-O(\exp(-n))$.

An upper bound follows similarly. By \cref{lem:mineig} we have that with probability at least $1-O(\exp(-n)$,
\begin{equation}
\begin{split}
n^{-1} x_j^{\top}(p^{-1}XX^{\top})^{-2}x_{j} \leq C\frac{p}{n}x_j^{\top}(XX^{\top})^{-1}x_j=C\frac{p}{n}e_{j}U_{n}U_n^{\top}e_{j}.
\end{split}
\end{equation}
Then, by \cref{lem:JL}, with probability at least $1-O(\exp(-n))$, $\hat{\Omega}_{jj}\leq C\frac{1}{1-\epsilon_{u}}$. By \cref{ass:sparsity} $n$ is sufficiently larger than $\log p$, so applying a union bound over $p$ gives $\max_{j=1,\ldots,p}\hat{\Omega}_{jj}<C$ for some positive constant $C$ with probability $1-O(\exp(-n))$. We conclude that $\hat{\Omega}_{jj}$ is $O_{p}(1)$ for $j=1,\ldots,p$.

\subsection{Proof of \cref{theorem:biasvanishridge}}\label{app:proofridgelemma}
Consider the $(i,j)$-th element of $\hat{M}X$,
\begin{equation}
    [\hat{M}X]_{i,j} = \left[e_{i}^{\top}(X^{\top}X+\gamma I_{p})^{-1}X^{\top}Xe_{i}\right]^{-1}e_{i}^{\top}(X^{\top}X + \gamma I_{p})^{-1}X^{\top}Xe_{j}.
\end{equation}
By substituting $X = VSU^{\top}=VS_{n}U_{n}^{\top}$, where $U^{\top}U=UU^{\top}=I_{p}$, $U_{n}^{\top}U_{n} =I_{p}$, ${Q} = S^{\top}S$, and $Q_{n}=S_{n}^{\top}S_{n}$ we have
\begin{equation}
\begin{split}
[\hat{M}X]_{i,j}& = \left[e_{i}^{\top}(U {Q} U^{\top}+\gamma I_{p})^{-1}U{Q} U^{\top}e_{i}\right]^{-1}e_{i}^{\top}(U {Q} U^{\top}+\gamma I_{p})^{-1}U{Q} U^{\top}e_{j}\\
& = \left[e_{i}^{\top}(U_{n}U_{n}^{\top} + U_{n}(W-I_{n})U_{n}^{\top})e_{i}\right]^{-1}e_{i}^{\top}(U_{n}U_{n}^{\top} + U_{n}(W-I_{n})U_{n}^{\top})e_{j}\\
&=\frac{e_{i}^{\top}U_{n}U_{n}^{\top}e_{j}}{e_{i}^{\top}U_{n}U_{n}^{\top}e_{i}}\frac{1}{1+\frac{e_{i}^{\top}U_{n}(W-I_{n})U_{n}^{\top}e_{i}}{e_{i}^{\top}U_{n}U_{n}^{\top}e_{i}}} + \frac{e_{i}^{\top}U_{n}(W-I_{n})U_{n}^{\top}e_{j}}{e_{i}^{\top}U_{n}U_{n}^{\top}e_{i}\left(1+\frac{e_{i}^{\top}U_{n}(W-I_{n})U_{n}^{\top}e_{i}}{e_{i}^{\top}U_{n}U_{n}^{\top}e_{i}}\right)},\notag
\end{split}
\end{equation}
where $U_{n}U_{n}^{\top}=X^{\top}(XX^{\top})^{-1}X^{\top}$, $[W]_{ii} = \frac{{q}_{i}}{{q}_{i}+\gamma}$, and $[W]_{ii}-1 = -\frac{\gamma}{\gamma+{q}_{i}}$, for $i=1,\ldots,n$.

Now, 
\begin{equation}
\begin{split}
\frac{e_{i}^{\top}U_{n}(W-I_{n})U_{n}^{\top}e_{i}}{e_{i}^{\top}U_{n}U_{n}^{\top}e_{i}} &\leq \max_{i=1,\ldots,n}\left|\frac{\gamma}{\gamma+{q}_{i}}\right|\frac{e_{i}^{\top}U_{n}U_{n}^{\top}e_{i}}{e_{i}^{\top}U_{n}U_{n}^{\top}e_{i}}\\
&= \left|\frac{\gamma}{\gamma+\lambda_{\min\neq 0}(X^{\top}X)}\right|\\
&\leq \frac{\gamma}{\gamma + Cp},
\end{split}
\end{equation}
with the last line holding with probability $1-O(\exp(-n))$ by \cref{lem:mineig}. Choosing $\gamma = O\left(p\sqrt{\frac{\log (p)}{n}}\right)$, we have that with probability $1-O(\exp(-n))$,
\begin{equation}
\begin{split}
\frac{e_{i}^{\top}U_{n}(W-I_{n})U_{n}^{\top}e_{i}}{e_{i}^{\top}U_{n}U_{n}^{\top}e_{i}} &=O(\sqrt{\log p/n})
\end{split}
\end{equation}
Then, by Cauchy-Schwarz, $e_{i}^{\top}U_{n}(W-I_{n})U_{n}^{\top}e_{j}/(e_{i}^{\top}U_{n}U_{n}^{\top}e_{i}) = O(\sqrt{\log p/n})$ with probability $1-O(\exp(-n))$. In that case,  $[MX]_{i,j} = \frac{e_{i}^{\top}U_{n}U_{n}^{\top}e_{j}}{e_{i}^{\top}U_{n}U_{n}^{\top}e_{i}} + O(\sqrt{\log p/n})$ with probability $1-O(\exp(-n))$, and \cref{theorem:biasvanishridge} then follows from \cref{theorem:biasvanish}.

\subsection{Proof of \cref{theorem:finitevarridge}}\label{app:prooftheoremfinitevarridge}
Normality follows from \cref{ass:error}. What remains is to show that $\hat{\Omega}_{jj}=O_{p}(1)$.
\begin{equation}\label{eq:varOmegaridge}
    \hat{\Omega}_{jj} =\text{var}(Z_{j}|X)=\sigma^2 \frac{\frac{1}{n}e_{j}^{\top}(X^{\top}X+\gamma I_{p})^{-1}X^{\top}X (X^{\top}X+\gamma I_{p})^{-1}e_{j}}{(\frac{1}{n}e_{j}^{\top}(X^{\top}X+\gamma I_{p})^{-1}X^{\top}Xe_{j})^{2}}.
\end{equation}
First the lower bound. By Cauchy-Schwarz
\begin{equation}
\begin{split}
    \left[\frac{1}{n}e_{j}^{\top}(X^{\top}X+\gamma I_{p})^{-1}X^{\top}Xe_{j}\right]^{2}&\leq \left[\frac{1}{n}e_{j}^{\top}X^{\top}Xe_{j}\right]\\
    &\quad \times\left[ \frac{1}{n}e_{j}^{\top}(X^{\top}X+\gamma I_{p})^{-1}X^{\top}X (X^{\top}X+\gamma I_{p})^{-1}e_{j}\right].
    \end{split}
\end{equation}
We therefore have
\begin{equation}
    \hat{\Omega}_{jj}\geq \left(\frac{1}{n}e_{j}^{\top}X^{\top}Xe_{j}\right)^{-1},
\end{equation}
and the lower bound follows from \cref{lem:use}.

We now continue with an upper bound.
Define $\bar{W}$ as a diagonal matrix with $[\bar{W}]_{ii}=\frac{1}{\gamma+{q}_{i}}$. Note that by \cref{lem:mineig}, $[\bar{W}]_{ii}\leq Cp^{-1}$ with probability $1-O(\exp(-n))$. The numerator of \eqref{eq:varOmegaridge} is then equal to 
\begin{equation}
\begin{split}
    \frac{1}{n}e_{j}^{\top}U_{n}\bar{W}(I_{n}-W)U_{n}^{\top}e_{j}&\leq C\frac{1}{pn}e_{j}^{\top}U_{n}(I_{n}-W)U_{n}^{\top}e_{j}\\
    &=C\frac{1}{pn}e_{j}^{\top}(X^{\top}X+\gamma I_{p})^{-1}X^{\top}Xe_{j},
    \end{split}
\end{equation}
with probability $1-O(\exp(-n))$.

From this we have that with probability   $1-O(\exp(-n))$,
\begin{equation}
    \hat{\Omega}_{jj}\leq \sigma^2 C\frac{1}{p} \frac{1}{\frac{1}{n}e_{j}^{\top}U_{n}(I_{n}-W)U_{n}^{\top}e_{j}}\leq \sigma^2 C_{1}\frac{1}{p} \frac{1}{\frac{1}{n}e_{j}^{\top}U_{n}U_{n}^{\top}e_{j}},
    \end{equation}
    where we used that $I_{n}-W = ({Q}_n + \gamma I_{n})^{-1}{Q}_n$ and 
    \begin{align*}
        \min_{i=1,\ldots,n}e_{i}^{\top}(I_{n}-W)e_{i}\geq \frac{\lambda_{\min\neq 0}(X^{\top}X)}{\lambda_{\min\neq 0}(X^{\top}X)+\gamma} \geq \frac{p}{p(1+\sqrt{\log p/n})}\geq c.
    \end{align*}
    with the last inequality following from \cref{ass:sparsity}. 
    Under \cref{ass:symmetry}, we can now apply \cref{lem:JL} and find that with probability $1-O(\exp(-n))$,
    \begin{equation}
        \hat{\Omega}_{jj}\leq \sigma^2C.
\end{equation}
The proof now follows from a union bound.

\clearpage
\section{Results Monte Carlo experiments }\label{A: mcresults}
\begin{table}[H]
  \centering \footnotesize
  \caption{Results independent and equicorrelated designs}
    \begin{tabular}{llrrrrrrrrr}
    \toprule \toprule
          &       & \multicolumn{3}{c}{(p,n)=(200,100)} & \multicolumn{3}{c}{(p,n)=(1000,200)} & \multicolumn{3}{c}{(p,n)=(10000,400)} \\
          &       & \multicolumn{1}{l}{MAE} & \multicolumn{1}{l}{CR} & \multicolumn{1}{l}{power} & \multicolumn{1}{l}{MAE} & \multicolumn{1}{l}{CR} & \multicolumn{1}{l}{power} & \multicolumn{1}{l}{MAE} & \multicolumn{1}{l}{CR} & \multicolumn{1}{l}{power} \\
         \cmidrule(lr){3-5}\cmidrule(lr){6-8}\cmidrule(lr){9-11}
method & $S$ & \multicolumn{9}{c}{Independent} \\
\midrule
    \multirow{2}[0]{*}{MPI} & $S$   & 0.70  & 0.98  & 0.63  & 0.38  & 0.96  & 0.97  & 0.22  & 0.95  & 1.00 \\
          & $S^c$ & 0.64  & 0.99  &       & 0.30  & 0.99  &       & 0.18  & 0.98  &  \\[+1mm]
    \multirow{2}[0]{*}{RID} & $S$   & 0.62  & 0.97  & 0.77  & 0.37  & 0.95  & 0.97  & 0.22  & 0.95  & 1.00 \\
          & $S^c$ & 0.52  & 0.99  &       & 0.29  & 0.99  &       & 0.17  & 0.98  &  \\[+1mm]
    \multirow{2}[0]{*}{GBRD} & $S$   & 0.58  & 0.96  & 0.85  & 0.36  & 0.94  & 0.98  &       &       &  \\
          & $S^c$ & 0.43  & 0.99  &       & 0.26  & 0.99  &       &       &       &  \\[+1mm]
    \multirow{2}[0]{*}{JM} & $S$   & 0.66  & 0.84  & 0.89  & 0.41  & 0.83  & 0.99  &       &       &  \\
          & $S^c$ & 0.33  & 0.98  &       & 0.21  & 0.99  &       &       &       &  
 \\\midrule
          &       & \multicolumn{9}{c}{Equicorrelated with $\rho=0.3$} \\ \midrule
\multirow{2}[0]{*}{MPI} & $S$   & 2.31  & 0.95  & 0.37  & 1.38  & 0.93  & 0.81  & 0.90  & 0.91  & 0.99 \\
          & $S^c$ & 2.11  & 0.97  &       & 1.15  & 0.97  &       & 0.74  & 0.96  &  \\[+1mm]
    \multirow{2}[0]{*}{RID} & $S$   & 1.94  & 0.93  & 0.51  & 1.34  & 0.92  & 0.84  & 0.90  & 0.91  & 0.99 \\
          & $S^c$ & 1.65  & 0.97  &       & 1.08  & 0.97  &       & 0.73  & 0.96  &  \\[+1mm]
    \multirow{2}[0]{*}{GBRD} & $S$   & 1.82  & 0.91  & 0.59  & 1.32  & 0.90  & 0.86  &       &       &  \\
          & $S^c$ & 1.37  & 0.97  &       & 0.98  & 0.97  &       &       &       &  \\[+1mm]
    \multirow{2}[0]{*}{JM} & $S$   & 2.00  & 0.59  & 0.72  & 1.71  & 0.51  & 0.91  &       &       &  \\
          & $S^c$ & 0.82  & 0.95  &       & 0.57  & 0.97  &       &       &       &  \\\midrule
          &       & \multicolumn{9}{c}{Equicorrelated with $\rho=0.6$} \\ \midrule
    \multirow{2}[0]{*}{MPI} & $S$   & 3.66  & 0.94  & 0.19  & 2.11  & 0.94  & 0.43  & 1.39  & 0.94  & 0.79 \\
          & $S^c$ & 3.51  & 0.95  &       & 1.90  & 0.96  &       & 1.23  & 0.96  &  \\[+1mm]
    \multirow{2}[0]{*}{RID} & $S$   & 2.86  & 0.94  & 0.26  & 2.00  & 0.94  & 0.47  & 1.38  & 0.94  & 0.80 \\
          & $S^c$ & 2.60  & 0.96  &       & 1.76  & 0.96  &       & 1.21  & 0.96  &  \\[+1mm]
    \multirow{2}[0]{*}{GBRD} & $S$   & 2.68  & 0.93  & 0.28  & 1.96  & 0.93  & 0.48  &       &       &  \\
          & $S^c$ & 2.25  & 0.96  &       & 1.61  & 0.96  &       &       &       &  \\[+1mm]
    \multirow{2}[0]{*}{JM} & $S$   & 2.99  & 0.31  & 0.39  & 2.80  & 0.25  & 0.51  &       &       &  \\
          & $S^c$ & 0.80  & 0.95  &       & 0.61  & 0.98  &       &       &       &  \\
\midrule
& & \multicolumn{9}{c}{Equicorrelated with $\rho=0.9$}\\
\midrule
    \multirow{2}[0]{*}{MPI} & $S$   & 8.37  & 0.95  & 0.08  & 4.60  & 0.96  & 0.12  & 2.90  & 0.96  & 0.23 \\
          & $S^c$ & 8.28  & 0.96  &       & 4.53  & 0.96  &       & 2.92  & 0.96  &  \\[+1mm]
    \multirow{2}[0]{*}{RID} & $S$   & 5.70  & 0.95  & 0.09  & 4.07  & 0.96  & 0.12  & 2.83  & 0.96  & 0.23 \\
          & $S^c$ & 5.56  & 0.96  &       & 3.99  & 0.96  &       & 2.85  & 0.96  &  \\[+1mm]
    \multirow{2}[0]{*}{GBRD} & $S$   & 5.49  & 0.95  & 0.09  & 3.96  & 0.95  & 0.12  &       &       &  \\
          & $S^c$ & 5.28  & 0.96  &       & 3.85  & 0.96  &       &       &       &  \\[+1mm]
    \multirow{2}[0]{*}{JM} & $S$   & 4.17  & 0.10  & 0.12  & 3.57  & 0.47  & 0.13  &       &       &  \\
          & $S^c$ & 0.80  & 0.96  &       & 1.28  & 0.97  &       &       &       &  \\
          \bottomrule \bottomrule
    \end{tabular}%
  \label{tab:independent_equicorrelated}%
\end{table}%
\begin{table}[H]
  \centering \footnotesize
  \caption{Results auto-regressive designs}
    \begin{tabular}{llrrrrrrrrr}
    \toprule \toprule
          &       & \multicolumn{3}{c}{(p,n)=(200,100)} & \multicolumn{3}{c}{(p,n)=(1000,200)} & \multicolumn{3}{c}{(p,n)=(10000,400)} \\
          &       & \multicolumn{1}{l}{MAE} & \multicolumn{1}{l}{CR} & \multicolumn{1}{l}{power} & \multicolumn{1}{l}{MAE} & \multicolumn{1}{l}{CR} & \multicolumn{1}{l}{power} & \multicolumn{1}{l}{MAE} & \multicolumn{1}{l}{CR} & \multicolumn{1}{l}{power} \\
         \cmidrule(lr){3-5}\cmidrule(lr){6-8}\cmidrule(lr){9-11}
method & $S$ & \multicolumn{9}{c}{Auto-regressive with $\rho=0.3$} \\
\midrule
    \multirow{2}[0]{*}{MPI} & $S$   & 0.48  & 0.97  & 0.81  & 0.27  & 0.96  & 1.00  & 0.17  & 0.94  & 1.00 \\
          & $S^c$ & 0.45  & 0.98  &       & 0.24  & 0.98  &       & 0.15  & 0.97  &  \\[+1mm]
    \multirow{2}[0]{*}{RID} & $S$   & 0.40  & 0.96  & 0.91  & 0.26  & 0.95  & 1.00  & 0.17  & 0.94  & 1.00 \\
          & $S^c$ & 0.36  & 0.98  &       & 0.23  & 0.98  &       & 0.15  & 0.97  &  \\[+1mm]
    \multirow{2}[0]{*}{GBRD} & $S$   & 0.37  & 0.95  & 0.95  & 0.25  & 0.96  & 1.00  &       &       &  \\
          & $S^c$ & 0.30  & 0.98  &       & 0.21  & 0.98  &       &       &       &  \\[+1mm]
    \multirow{2}[0]{*}{JM} & $S$   & 0.41  & 0.80  & 0.96  & 0.30  & 0.78  & 1.00  &       &       &  \\
          & $S^c$ & 0.21  & 0.97  &       & 0.15  & 0.98  &       &       &       & \\
    \midrule
          &       & \multicolumn{9}{c}{Auto-regressive with $\rho=0.6$} \\ \midrule
    \multirow{2}[0]{*}{MPI} & $S$   & 0.64  & 0.97  & 0.65  & 0.34  & 0.95  & 0.98  & 0.21  & 0.91  & 1.00 \\
          & $S^c$ & 0.61  & 0.98  &       & 0.29  & 0.97  &       & 0.17  & 0.97  &  \\[+1mm]
    \multirow{2}[0]{*}{RID} & $S$   & 0.49  & 0.95  & 0.83  & 0.31  & 0.95  & 0.99  & 0.21  & 0.91  & 1.00 \\
          & $S^c$ & 0.43  & 0.98  &       & 0.27  & 0.97  &       & 0.17  & 0.97  &  \\[+1mm]
    \multirow{2}[0]{*}{GBRD} & $S$   & 0.47  & 0.95  & 0.86  & 0.31  & 0.95  & 0.98  &       &       &  \\
          & $S^c$ & 0.39  & 0.98  &       & 0.28  & 0.97  &       &       &       &  \\[+1mm]
    \multirow{2}[0]{*}{JM} & $S$   & 0.48  & 0.67  & 0.96  & 0.33  & 0.67  & 1.00  &       &       &  \\
          & $S^c$ & 0.20  & 0.97  &       & 0.15  & 0.97  &       &       &       &  \\
    \midrule
          &       & \multicolumn{9}{c}{Auto-regressive with $\rho=0.9$} \\ \midrule
    \multirow{2}[0]{*}{MPI} & $S$   & 1.60  & 0.95  & 0.23  & 0.78  & 0.90  & 0.63  & 0.49  & 0.74  & 0.98 \\
          & $S^c$ & 1.60  & 0.96  &       & 0.65  & 0.96  &       & 0.27  & 0.96  &  \\[+1mm]
    \multirow{2}[0]{*}{RID} & $S$   & 0.99  & 0.89  & 0.49  & 0.68  & 0.83  & 0.79  & 0.48  & 0.72  & 0.99 \\
          & $S^c$ & 0.78  & 0.96  &       & 0.47  & 0.96  &       & 0.25  & 0.96  &  \\[+1mm]
    \multirow{2}[0]{*}{GBRD} & $S$   & 1.08  & 0.91  & 0.46  & 0.76  & 0.91  & 0.66  &       &       &  \\
          & $S^c$ & 0.89  & 0.97  &       & 0.64  & 0.97  &       &       &       &  \\[+1mm]
    \multirow{2}[0]{*}{JM} & $S$   & 0.89  & 0.38  & 0.75  & 0.69  & 0.57  & 0.85  &       &       &  \\
          & $S^c$ & 0.23  & 0.97  &       & 0.27  & 0.96  &       &       &       &  \\
          \bottomrule \bottomrule
    \end{tabular}%
  \label{tab:autoregressive}%
\end{table}%
\begin{table}[H]
  \centering \footnotesize
  \caption{Results factor model designs}
    \begin{tabular}{llrrrrrrrrr}
    \toprule \toprule
          &       & \multicolumn{3}{c}{(p,n)=(200,100)} & \multicolumn{3}{c}{(p,n)=(1000,200)} & \multicolumn{3}{c}{(p,n)=(10000,400)} \\
          &       & \multicolumn{1}{l}{MAE} & \multicolumn{1}{l}{CR} & \multicolumn{1}{l}{power} & \multicolumn{1}{l}{MAE} & \multicolumn{1}{l}{CR} & \multicolumn{1}{l}{power} & \multicolumn{1}{l}{MAE} & \multicolumn{1}{l}{CR} & \multicolumn{1}{l}{power} \\
         \cmidrule(lr){3-5}\cmidrule(lr){6-8}\cmidrule(lr){9-11}
method & $S$ & \multicolumn{9}{c}{Factor model with $k=2$} \\
\midrule
    \multirow{2}[0]{*}{MPI} & $S$   & 2.14  & 0.97  & 0.40  & 1.25  & 0.95  & 0.79  & 0.88  & 0.93  & 0.95 \\
          & $S^c$ & 2.00  & 0.98  &       & 1.09  & 0.98  &       & 0.67  & 0.98  &  \\[+1mm]
    \multirow{2}[0]{*}{RID} & $S$   & 1.70  & 0.96  & 0.58  & 1.19  & 0.95  & 0.83  & 0.87  & 0.93  & 0.96 \\
          & $S^c$ & 1.44  & 0.98  &       & 0.99  & 0.98  &       & 0.66  & 0.98  &  \\[+1mm]
    \multirow{2}[0]{*}{GBRD} & $S$   & 1.70  & 0.91  & 0.63  & 1.29  & 0.87  & 0.84  &       &       &  \\
          & $S^c$ & 1.13  & 0.98  &       & 0.81  & 0.98  &       &       &       &  \\[+1mm]
    \multirow{2}[0]{*}{JM} & $S$   & 2.47  & 0.34  & 0.68  & 2.09  & 0.30  & 0.85  &       &       &  \\
          & $S^c$ & 0.45  & 0.96  &       & 0.37  & 0.98  &       &       &       &  \\
    \midrule
          &       & \multicolumn{9}{c}{Factor model with $k=10$} \\ \midrule
    \multirow{2}[0]{*}{MPI} & $S$   & 3.87  & 0.97  & 0.13  & 2.26  & 0.96  & 0.37  & 1.47  & 0.96  & 0.68 \\
          & $S^c$ & 3.76  & 0.98  &       & 2.10  & 0.97  &       & 1.37  & 0.97  &  \\[+1mm]
    \multirow{2}[0]{*}{RID} & $S$   & 2.57  & 0.94  & 0.27  & 1.98  & 0.94  & 0.44  & 1.44  & 0.95  & 0.70 \\
          & $S^c$ & 2.06  & 0.98  &       & 1.66  & 0.97  &       & 1.28  & 0.97  &  \\[+1mm]
    \multirow{2}[0]{*}{GBRD} & $S$   & 2.66  & 0.92  & 0.27  & 2.09  & 0.88  & 0.44  &       &       &  \\
          & $S^c$ & 1.99  & 0.98  &       & 1.45  & 0.97  &       &       &       &  \\[+1mm]
    \multirow{2}[0]{*}{JM} & $S$   & 3.40  & 0.15  & 0.43  & 3.35  & 0.10  & 0.50  &       &       &  \\
          & $S^c$ & 0.36  & 0.96  &       & 0.33  & 0.98  &       &       &       &  \\
    \midrule
          &       & \multicolumn{9}{c}{Factor model with $k=20$} \\ \midrule
    \multirow{2}[0]{*}{MPI} & $S$   & 5.08  & 0.98  & 0.06  & 2.97  & 0.97  & 0.18  & 1.99  & 0.96  & 0.43 \\
          & $S^c$ & 4.96  & 0.98  &       & 2.89  & 0.97  &       & 1.89  & 0.97  &  \\[+1mm]
    \multirow{2}[0]{*}{RID} & $S$   & 2.68  & 0.91  & 0.23  & 2.41  & 0.92  & 0.28  & 1.89  & 0.93  & 0.47 \\
          & $S^c$ & 1.85  & 0.98  &       & 1.84  & 0.97  &       & 1.59  & 0.97  &  \\[+1mm]
    \multirow{2}[0]{*}{GBRD} & $S$   & 3.00  & 0.95  & 0.16  & 2.52  & 0.91  & 0.26  &       &       &  \\
          & $S^c$ & 2.43  & 0.98  &       & 1.88  & 0.98  &       &       &       &  \\[+1mm]
    \multirow{2}[0]{*}{JM} & $S$   & 3.13  & 0.22  & 0.49  & 3.61  & 0.07  & 0.44  &       &       &  \\
          & $S^c$ & 0.38  & 0.96  &       & 0.20  & 0.99  &       &       &       &  \\
          \bottomrule \bottomrule
    \end{tabular}%
  \label{tab:factormodel}%
\end{table}%
\begin{table}[H]
  \centering \footnotesize
  \caption{Results group structure and extreme designs}
    \begin{tabular}{llrrrrrrrrr}
    \toprule \toprule
          &       & \multicolumn{3}{c}{(p,n)=(200,100)} & \multicolumn{3}{c}{(p,n)=(1000,200)} & \multicolumn{3}{c}{(p,n)=(10000,400)} \\
          &       & \multicolumn{1}{l}{MAE} & \multicolumn{1}{l}{CR} & \multicolumn{1}{l}{power} & \multicolumn{1}{l}{MAE} & \multicolumn{1}{l}{CR} & \multicolumn{1}{l}{power} & \multicolumn{1}{l}{MAE} & \multicolumn{1}{l}{CR} & \multicolumn{1}{l}{power} \\
         \cmidrule(lr){3-5}\cmidrule(lr){6-8}\cmidrule(lr){9-11}
method & $S$ & \multicolumn{9}{c}{Group structure with $\delta^2=0.01$} \\
\midrule
    \multirow{2}[0]{*}{MPI} & $S$   & 4.56  & 0.85  & 0.23  & 3.28  & 0.76  & 0.51  & 2.92  & 0.70  & 0.92 \\
          & $S^c$ & 2.85  & 0.98  &       & 1.50  & 0.98  &       & 0.98  & 0.97  &  \\[+1mm]
    \multirow{2}[0]{*}{RID} & $S$   & 3.78  & 0.81  & 0.27  & 3.22  & 0.75  & 0.54  & 2.91  & 0.70  & 0.92 \\
          & $S^c$ & 2.24  & 0.98  &       & 1.43  & 0.98  &       & 0.97  & 0.97  &  \\[+1mm]
    \multirow{2}[0]{*}{GBRD} & $S$   & 6.20  & 0.90  & 0.14  & 5.65  & 0.88  & 0.18  &       &       &  \\
          & $S^c$ & 1.86  & 0.98  &       & 1.34  & 0.98  &       &       &       &  \\[+1mm]
    \multirow{2}[0]{*}{JM} & $S$   & 3.00  & 0.79  & 0.27  & 3.72  & 0.79  & 0.26  &       &       &  \\
          & $S^c$ & 1.44  & 0.98  &       & 1.07  & 0.98  &       &       &       &  \\
   \midrule
          &       & \multicolumn{9}{c}{Group structure with $\delta^2=0.05$} \\ 
          \midrule
           \multirow{2}[0]{*}{MPI} & $S$   & 4.73  & 0.90  & 0.18  & 3.14  & 0.73  & 0.52  & 2.80  & 0.62  & 0.92 \\
          & $S^c$ & 2.76  & 0.98  &       & 1.46  & 0.98  &       & 0.96  & 0.97  &  \\[+1mm]
    \multirow{2}[0]{*}{RID} & $S$   & 3.75  & 0.84  & 0.27  & 3.08  & 0.72  & 0.55  & 2.80  & 0.62  & 0.92 \\
          & $S^c$ & 2.18  & 0.98  &       & 1.40  & 0.98  &       & 0.95  & 0.97  &  \\[+1mm]
    \multirow{2}[0]{*}{GBRD} & $S$   & 5.11  & 0.94  & 0.12  & 4.16  & 0.90  & 0.19  &       &       &  \\
          & $S^c$ & 1.81  & 0.98  &       & 1.31  & 0.98  &       &       &       &  \\[+1mm]
    \multirow{2}[0]{*}{JM} & $S$   & 2.75  & 0.77  & 0.33  & 3.15  & 0.79  & 0.33  &       &       &  \\
          & $S^c$ & 1.41  & 0.98  &       & 1.05  & 0.98  &       &       &       &  \\
\midrule
          &       & \multicolumn{9}{c}{Group structure with $\delta^2=0.1$} \\ 
          \midrule
    \multirow{2}[0]{*}{MPI} & $S$   & 4.55  & 0.92  & 0.16  & 2.99  & 0.74  & 0.53  & 2.69  & 0.59  & 0.92 \\
          & $S^c$ & 2.68  & 0.98  &       & 1.43  & 0.98  &       & 0.94  & 0.97  &  \\[+1mm]
    \multirow{2}[0]{*}{RID} & $S$   & 3.63  & 0.86  & 0.26  & 2.95  & 0.73  & 0.57  & 2.69  & 0.59  & 0.93 \\
          & $S^c$ & 2.13  & 0.98  &       & 1.37  & 0.98  &       & 0.93  & 0.97  &  \\[+1mm]
    \multirow{2}[0]{*}{GBRD} & $S$   & 4.28  & 0.94  & 0.13  & 3.42  & 0.91  & 0.21  &       &       &  \\
          & $S^c$ & 1.77  & 0.98  &       & 1.28  & 0.98  &       &       &       &  \\[+1mm]
    \multirow{2}[0]{*}{JM} & $S$   & 2.51  & 0.77  & 0.39  & 2.64  & 0.75  & 0.43  &       &       &  \\
          & $S^c$ & 1.37  & 0.98  &       & 1.03  & 0.98  &       &       &       &  \\
\midrule
& & \multicolumn{9}{c}{Extreme correlation}\\
\midrule
    \multirow{2}[0]{*}{MPI} & $S$   & 3.05  & 0.85  & 0.19  & 2.66  & 0.36  & 0.41  & 2.62  & 0.33  & 0.53 \\
          & $S^c$ & 3.06  & 0.98  &       & 1.63  & 0.98  &       & 1.06  & 0.98  &  \\[+1mm]
    \multirow{2}[0]{*}{RID} & $S$   & 2.82  & 0.42  & 0.31  & 2.76  & 0.33  & 0.40  & 2.69  & 0.31  & 0.52 \\
          & $S^c$ & 2.09  & 0.97  &       & 1.46  & 0.98  &       & 1.04  & 0.98  &  \\[+1mm]
    \multirow{2}[0]{*}{GBRD} & $S$   & 4.03  & 0.69  & 0.16  & 3.75  & 0.75  & 0.19  &       &       &  \\
          & $S^c$ & 1.91  & 0.98  &       & 1.36  & 0.98  &       &       &       &  \\[+1mm]
    \multirow{2}[0]{*}{JM} & $S$   & 3.62  & 0.17  & 0.34  & 3.48  & 0.29  & 0.37  &       &       &  \\
          & $S^c$ & 0.35  & 0.95  &       & 0.48  & 0.98  &       &       &       &  \\
          \bottomrule \bottomrule
    \end{tabular}%
  \label{tab:group_extreme}%
\end{table}%

\end{appendices}

\end{document}